\documentclass[10pt,a4paper]{article}
\usepackage{CJK}
\usepackage{mathrsfs}
\usepackage{amsfonts}
\begin{document}
\begin{CJK*}{GBK}{song}
\begin{center}
 {\huge\textbf{New Properties of Fourier Series and Riemann Zeta Function}}
\vspace{0.3cm}

\footnotetext{\hspace*{-.45cm}\footnotesize $^\dag$ Corresponding
author. E-mail: yananbiguangqing@sohu.com }

\begin{center}
\rm Guangqing Bi $^{\rm a)\dagger}$, \ \ Yuekai Bi $^{\rm b)}$
\end{center}

\begin{center}
\begin{footnotesize} \sl
${}^{\rm b)}$ School of Electronic and Information Engineering,
BUAA, Beijing 100191, China\\
E-mail: yuekaifly@163.com
\end{footnotesize}
\end{center}

\end{center}

\noindent
\begin{abstract}
We establish the mapping relations between analytic functions and
periodic functions using the abstract operators $\cos(h\partial_x)$
and $\sin(h\partial_x)$, including the mapping relations between
power series and trigonometric series, and by using such mapping
relations we obtain a general method to find the sum function of a
trigonometric series. According to this method, if each coefficient
of a power series is respectively equal to that of a trigonometric
series, then if we know the sum function of the power series, we can
obtain that of the trigonometric series, and the non-analytical
points of which are also determined at the same time, thus we obtain
a general method to find the sum of the Dirichlet series of integer
variables, and derive several new properties of $\zeta(2n+1)$.
\end{abstract}


\noindent{\bf MSC(2000) Subject Classification}\ \ 11M06; 42A24;
11M35; 35S99

\section{Introduction}
\noindent

The trigonometric series, especially the Fourier series, is of great
importance in both mathematics and physics, and expanding periodic
functions into Fourier series has become a very mature theory.
During the early years of last century, people realized the
significance of the inverse problem, which is how to find the sum
function of a certain Fourier series. So the question is: can we use
the sum function of power series to obtain that of trigonometric
series? As we know, the domain of functions expressed by
trigonometric series can be extended into the entire number axis
with the existence of denumerable non-analytical points, thus
functions expressed by trigonometric series are piecewise analytic
periodic functions. If a trigonometric series converges to an
analytic sum function in a certain interval, then it is quite
natural that there is a mapping relation between a periodic function
and an analytic function, though the two endpoints of the interval
generally are its non-analytical points. Therefore, we can use the
sum of power series to obtain that of corresponding trigonometric
series, converting the research focus from trigonometric series to
power series. For a long time, the first author has realized the
fact that the abstract operators $\cos(h\partial_x)$ and
$\sin(h\partial_x)$ can express the mapping relations between
periodic functions and analytic functions more distinctly, which has
become a significant tool to obtain the sum function of
trigonometric series. Why the abstract operators $\cos(h\partial_x)$
and $\sin(h\partial_x)$ can establish such relations? That is
because this kind of operators is also a kind of functions,
containing the duality between periodic functions and linear
operators. By using such a duality, the authors have preliminarily
established theories of partial differential equations of abstract
operators in reference \cite{bi97}-\cite{bi10'}.

What is an abstract operator? The operator $f(t,\partial_x)$ is
generally interpreted as a Taylor expansion, called the infinite
order operator. However, the need to consider its disk of
convergence is a major constraint of its broad applications.
Therefore, the first author has defined the operator
$f(t,\partial_x)$ as $f(t,\partial_x)e^{ax}=f(t,a)e^{ax}$ in 1997.
Each operator $f(t,\partial_x)$ has a set of algorithms without the
need to use its Taylor expansion, and the first author has also
provided a method to determine such kind of algorithms. In this
sense, the operator $f(t,\partial_x)$ is known as the abstract
operator in reference \cite{bi97}. However, the concept of abstract
operators has not yet been spread adequately. As a result, it is
easy to mistake the abstract operator for the infinite order
operator by the similar symbol they are using. In fact, only several
simple abstract operators can be expanded into Taylor series under
certain conditions, such as $\cos(h\partial_x)$ and
$\sin(h\partial_x)$. However, the particular method we use to
calculate is not by using their Taylor expansions but the certain
algorithms of $\cos(h\partial_x)$ and $\sin(h\partial_x)$
established by the authors in their former published papers, such as
the formulas (\ref{3}), (\ref{6}), (\ref{17}), (\ref{19}) and
(\ref{20}). Although there may be several non-analytical points in
the process of calculation, making the results at these points not
tenable, yet we will not encounter the infinite series, which
benefits us a lot by avoiding the need to consider the disk of
convergence.

The abstract operators $\cos(h\partial_x)$ and $\sin(h\partial_x)$
have exclusive advantages in symbol expressions, and the summation
method of Fourier series established from which can be directly
extended into multiple Fourier series.

By using this summation method of Fourier series, the abstract
operators can also have significant applications in the Riemann Zeta
function $\zeta(m)$ of the analytic number theory. The Riemann Zeta
function $\zeta(s)$ defined usually by the Dirichlet series
\begin{equation}\label{00}
\zeta(s)=\sum^\infty_{n=1}\frac{1}{n^s},\qquad(\Re(s)>1).
\end{equation}

In 1735, Euler proved that for an arbitrary even number $2K>0$,
$\zeta(2K)=a_{2K}\pi^{2K}$, where $a_{2K}$ is a rational number.
However, for all odd numbers $2K+1$, the arithmetic properties of
$\zeta(2K+1)$ are still unknown. As the summation method of
trigonometric series of abstract operators is quite suitable to find
the sum of the Dirichlet series of integer variables including the
Zeta function, we take an important step in studying the arithmetic
properties of $\zeta(2K+1)$.

\section{Basic formulas of abstract operators}
\noindent

When acting on elementary functions, the abstract operators
$\cos(h\partial_x)$ and $\sin(h\partial_x)$ have complete basic
formulas as differential operations, now we are going to use the
algorithms of abstract operators in reference \cite{bi97} or
\cite{bi10} to establish these formulas.

Firstly, according to the definition of abstract operators, we have
\begin{equation}\label{1}
    \cos(h\partial_x)e^{bx}=\cos(bh)e^{bx},\qquad\sin(h\partial_x)e^{bx}=\sin(bh)e^{bx}.
\end{equation}
Where
$bx=b_1x_1+b_2x_2+\cdots+b_nx_n,\;bh=b_1h_1+b_2h_2+\cdots+b_nh_n$.

According to Theorem 2 in reference \cite{bi97}, namely
\[\exp(ih\partial_x)f(x)=f(x+ih),\] and Theorem 3:
\[\cos(h\partial_x)f(x)=\Re[f(x+ih)],\qquad\sin(h\partial_x)f(x)=\Im[f(x+ih)],\]
we have

\parbox{11cm}{\begin{eqnarray*}\label{2}
                \cos(h\partial_x)\cos{bx} &=& \cosh(bh)\cos{bx}, \\
                \cos(h\partial_x)\sin{bx} &=& \cosh(bh)\sin{bx}, \\
                \sin(h\partial_x)\cos{bx} &=& -\sinh(bh)\sin{bx}, \\
                \sin(h\partial_x)\sin{bx} &=& \sinh(bh)\cos{bx}.
              \end{eqnarray*}}\hfill\parbox{1cm}{\begin{eqnarray}\end{eqnarray}}

Based on (\ref{2}), and by using Theorem 6 in reference \cite{bi97}:

\parbox{11cm}{\begin{eqnarray*}\label{3}
                \sin(h\partial_x)\frac{u}{v} &=& \frac{\cos(h\partial_x)v\cdot\sin(h\partial_x)u-\sin(h\partial_x)v\cdot\cos(h\partial_x)u}
                {(\cos(h\partial_x)v)^2+(\sin(h\partial_x)v)^2}, \\
                \cos(h\partial_x)\frac{u}{v} &=& \frac{\cos(h\partial_x)v\cdot\cos(h\partial_x)u+\sin(h\partial_x)v\cdot\sin(h\partial_x)u}
                {(\cos(h\partial_x)v)^2+(\sin(h\partial_x)v)^2},
              \end{eqnarray*}}\hfill\parbox{1cm}{\begin{eqnarray}\end{eqnarray}}
we have

\parbox{10cm}{\begin{eqnarray*}\label{4}
                \cos(h\partial_x)\tan{bx} &=& \frac{\sin(2bx)}{\cosh(2bh)+\cos(2bx)}, \\
                \cos(h\partial_x)\cot{bx} &=& \frac{\sin(2bx)}{\cosh(2bh)-\cos(2bx)},\\
                \sin(h\partial_x)\tan{bx} &=& \frac{\sinh(2bh)}{\cosh(2bh)+\cos(2bx)}, \\
                \sin(h\partial_x)\cot{bx} &=&
                \frac{\sinh(2bh)}{\cos(2bx)-\cosh(2bh)}.
              \end{eqnarray*}}\hfill\parbox{1cm}{\begin{eqnarray}\end{eqnarray}}

For secant and cosecant functions, similarly to (\ref{4}), we have

\parbox{10cm}{\begin{eqnarray*}\label{5}
                \cos(h\partial_x)\sec{bx} &=& \frac{2\cosh(bh)\cos{bx}}{\cosh(2bh)+\cos(2bx)}, \\
                \cos(h\partial_x)\csc{bx} &=& \frac{2\cosh(bh)\sin{bx}}{\cosh(2bh)-\cos(2bx)}, \\
                \sin(h\partial_x)\sec{bx} &=& \frac{2\sinh(bh)\sin{bx}}{\cosh(2bh)+\cos(2bx)}, \\
                \sin(h\partial_x)\csc{bx} &=&
                \frac{2\sinh(bh)\cos{bx}}{\cos(2bx)-\cosh(2bh)}.
              \end{eqnarray*}}\hfill\parbox{1cm}{\begin{eqnarray}\end{eqnarray}}

According to the following theorem:

\textbf{Theorem 1.} \cite{bi10} If $y=f(bx)\in{J}$ (set of analytic
functions) is the inverse function of $bx=g(y)$, namely
$g(f(bx))=bx$, then $\sin(h\partial_x)f(bx)$(denoted by $Y$) and
$\cos(h\partial_x)f(bx)$(denoted by $X$) can be determined by the
following set of equations:

\parbox{6cm}{\begin{eqnarray*}\label{6}
                \cos\left(Y\frac{\partial}{\partial{X}}\right)g(X) &=& bx \\
                \sin\left(Y\frac{\partial}{\partial{X}}\right)g(X) &=& bh
              \end{eqnarray*}}
              \hfill\parbox{4cm}{\begin{eqnarray*}x\in\mathbb{R}^n,\quad{h}\in\mathbb{R}_n,\end{eqnarray*}}
              \hfill\parbox{1cm}{\begin{eqnarray}\end{eqnarray}}

we can derive basic formulas of the corresponding inverse function:

\parbox{10cm}{\begin{eqnarray*}\label{7}
                \cos(h\partial_x)\ln(bx) &=& \ln\left((bx)^2+(bh)^2\right)^{1/2}, \\
                \sin(h\partial_x)\ln(bx) &=&\textrm{
                arccot}\frac{bx}{bh}.
              \end{eqnarray*}}\hfill\parbox{1cm}{\begin{eqnarray}\end{eqnarray}}

\parbox{10cm}{\begin{eqnarray*}\label{8}
                \sin(h\partial_x)\arctan{bx} &=& \frac{1}{2}\textrm{tanh}^{-1}\frac{2bh}{1+(bx)^2+(bh)^2}, \\
                \cos(h\partial_x)\arctan{bx} &=&
                \frac{1}{2}\arctan\frac{2bx}{1-(bx)^2-(bh)^2}.
              \end{eqnarray*}}\hfill\parbox{1cm}{\begin{eqnarray}\end{eqnarray}}

\parbox{10cm}{\begin{eqnarray*}\label{9}
                \sin(h\partial_x)\textrm{arccot}\,bx &=& -\frac{1}{2}\textrm{coth}^{-1}\frac{1+(bx)^2+(bh)^2}{2bh}, \\
                \cos(h\partial_x)\textrm{arccot}\,bx &=&
                \frac{1}{2}\textrm{arccot}\frac{(bx)^2+(bh)^2-1}{2bx}.
              \end{eqnarray*}}\hfill\parbox{1cm}{\begin{eqnarray}\end{eqnarray}}

\textbf{Proof.}  Here we only give the detailed proof of (\ref{7})
and (\ref{8}). According to (\ref{1}) and (\ref{6}), we have

\parbox{5cm}{\begin{eqnarray*}
                e^X\cos{Y} &=& bx \\
                e^X\sin{Y} &=& bh
              \end{eqnarray*}}\hfill\parbox{3cm}{\begin{eqnarray*} X &=& \cos(h\partial_x)\ln(bx), \\
                Y &=& \sin(h\partial_x)\ln(bx).\end{eqnarray*}}\hfill\parbox{4cm}{\begin{eqnarray*}\end{eqnarray*}}

By solving this set of equations we have (\ref{7}), and according to
(\ref{4}) and (\ref{6}), we have

\parbox{6cm}{\begin{eqnarray*}
                \frac{\sin2X}{\cosh2Y+\cos2X} &=& bx \\
                \frac{\sinh2Y}{\cosh2Y+\cos2X} &=& bh
              \end{eqnarray*}}\hfill\parbox{3cm}{\begin{eqnarray*} X &=& \cos(h\partial_x)\arctan{bx}, \\
                Y &=& \sin(h\partial_x)\arctan{bx}.\end{eqnarray*}}\hfill\parbox{3cm}{\begin{eqnarray*}\end{eqnarray*}}

By solving this set of equations we have
\[1+(bx)^2+(bh)^2=1+\frac{\sin^22X+\sinh^22Y}{(\cosh2Y+\cos2X)^2}=\frac{2\cosh2Y}{\cosh2Y+\cos2X}.\]
According to the second expression of this set of equations, we have
$\cosh2Y+\cos2X=(\sinh2Y)/bh$, and by substituting it into the above
expression, we have
\[1+(bx)^2+(bh)^2=\frac{2bh\cosh2Y}{\sinh2Y}\quad\mbox{or}\quad\tanh2Y=\frac{2bh}{1+(bx)^2+(bh)^2}.\]
Thus the first expression of (\ref{8}) is proved, similarly we can
prove the second one.

For hyperbolic and inverse hyperbolic functions, we can also derive
the corresponding basic formulas. For instance, correspondingly to
(\ref{2}), we have

\parbox{10cm}{\begin{eqnarray*}\label{10}
                \cos(h\partial_x)\cosh{bx} &=& \cos(bh)\cosh{bx}, \\
                \cos(h\partial_x)\sinh{bx} &=& \cos(bh)\sinh{bx}, \\
                \sin(h\partial_x)\cosh{bx} &=& \sin(bh)\sinh{bx}, \\
                \sin(h\partial_x)\sinh{bx} &=& \sin(bh)\cosh{bx}.
              \end{eqnarray*}}\hfill\parbox{1cm}{\begin{eqnarray}\end{eqnarray}}

Irrational functions can be considered as the inverse functions of
rational functions. For instance, $(bx)^{1/2}$ is the inverse
function of $(bx)^2$, thus similarly we have

\parbox{10cm}{\begin{eqnarray*}\label{11}
                \sin(h\partial_x)(bx)^{1/2} &=& \sqrt{\frac{\sqrt{(bx)^2+(bh)^2}-bx}{2}}, \\
                \cos(h\partial_x)(bx)^{1/2}
                &=&\sqrt{\frac{\sqrt{(bx)^2+(bh)^2}+bx}{2}}.
              \end{eqnarray*}}\hfill\parbox{1cm}{\begin{eqnarray}\end{eqnarray}}

To calculate concisely, the algorithms of products and composite
functions in reference \cite{bi97} are listed as follows, while the
algorithm of quotients has already been given in (\ref{3}):

\parbox{11cm}{\begin{eqnarray*}\label{19}
                \sin(h\partial_x)(vu) &=& \cos(h\partial_x)v\cdot\sin(h\partial_x)u+\sin(h\partial_x)v\cdot\cos(h\partial_x)u, \\
                \cos(h\partial_x)(vu) &=&
                \cos(h\partial_x)v\cdot\cos(h\partial_x)u-\sin(h\partial_x)v\cdot\sin(h\partial_x)u.
              \end{eqnarray*}}\hfill\parbox{1cm}{\begin{eqnarray}\end{eqnarray}}

\parbox{11cm}{\begin{eqnarray*}\label{20}
                \cos(h\partial_x)f(g(x)) &=& \cos\left(Y\frac{\partial}{\partial{X}} \right)f(X),\\
                \sin(h\partial_x)f(g(x)) &=& \sin\left(Y\frac{\partial}{\partial{X}}
                \right)f(X).
              \end{eqnarray*}}\hfill\parbox{1cm}{\begin{eqnarray}\end{eqnarray}}
Where $x\in\mathbb{R}^n$, $h\in\mathbb{R}_n$,
$X=\cos(h\partial_x)g(x)$, $Y=\sin(h\partial_x)g(x)$.

\section{Summation method of trigonometric series}
\noindent

\textbf{Theorem 2.}  Let $S(t)\in{J}$ be the sum function of the
power series $\sum^\infty_{n=0}a_nt^n$, $f(x)\in{L^2}[a,b]$ be the
sum function of the corresponding cosine series, and
$g(x)\in{L^2}[a,b]$ be that of the corresponding sine series, namely
\begin{eqnarray*}
  S(t) &=& \sum^\infty_{n=0}a_nt^n,\qquad{t}\in\mathbb{R}^1,\;\;0\leq{t}\leq{r},\;\;0<r<+\infty. \\
  f(x) &=& \sum^\infty_{n=0}a_n\cos\frac{n\pi{x}}{c}, \\
  g(x) &=&
  \sum^\infty_{n=0}a_n\sin\frac{n\pi{x}}{c},\qquad{x}\in\mathbb{R}^1,\;\;a<x<b,
\end{eqnarray*}
then we have the following mapping relations:

\parbox{11cm}{\begin{eqnarray*}\label{12}
                f(x) &=& \left.\cos\left(\frac{\pi{x}}{c}\frac{\partial}{\partial{z}}\right)S(e^z)\right|_{z=0}, \\
                g(x) &=&
                \left.\sin\left(\frac{\pi{x}}{c}\frac{\partial}{\partial{z}}\right)S(e^z)\right|_{z=0},
                \quad{z}\in\mathbb{R}^1.
              \end{eqnarray*}}\hfill\parbox{1cm}{\begin{eqnarray}\end{eqnarray}}
And the endpoints $a$ and $b$ of the interval $a<x<b$ are
non-analytical points (singularities) of Fourier series, which can
be uniquely determined by the detailed computation of the right-hand
side of (\ref{12}).

\textbf{Proof.}  By substituting $S(e^z)=\sum^\infty_{n=0}a_ne^{nz}$
into (\ref{12}), we can prove Theorem 2.

The sum function of infinite power series can be an elementary
function, which in most cases can be expressed as the definite
integral of an elementary function, thus in the application of
Theorem 2, the following theorem can be particularly useful:

\textbf{Theorem 3.}  Let $S(x)\in{J}$ be an arbitrary analytic
function integrable in the interval $[0,1]$, then we have

\parbox{10cm}{\begin{eqnarray*}\label{13}
 & & \left.\cos\left(\frac{\pi{x}}{c}\frac{\partial}{\partial{z}}\right)\int^{e^z}_0\!\!\!S(e^z)\,de^z\right|_{z=0}\\
 &=& \int^1_0\!\!S(\xi)\,d\xi-\frac{\pi}{c}\int^x_0\!\!\left.\sin\left(\frac{\pi{x}}{c}\frac{\partial}{\partial{z}}\right)
 [S(e^z)\,e^z]\right|_{z=0}dx. \\
\end{eqnarray*}}\hfill\parbox{1cm}{\begin{eqnarray}\end{eqnarray}}
\parbox{12cm}{\begin{eqnarray*}\label{14}
\left.\sin\left(\frac{\pi{x}}{c}\frac{\partial}{\partial{z}}\right)\int^{e^z}_0\!\!S(e^z)\,de^z\right|_{z=0}=
\frac{\pi}{c}\int^x_0\!\!\left.\cos\left(\frac{\pi{x}}{c}\frac{\partial}{\partial{z}}\right)[S(e^z)\,e^z]\right|_{z=0}dx.
\end{eqnarray*}}\hfill\parbox{1cm}{\begin{eqnarray}\end{eqnarray}}

\textbf{Proof.}  (\ref{13}) and (\ref{14}) are operator formulas.
According to the analytic continuous fundamental theorem in
reference \cite{bi97}, we only need to prove this set of formulas
when $S(x)=x^n$, $n\in\mathbb{N}_0$, this is obvious.

\textbf{Theorem 4.}   Let $S(x)\in{J}$ be an arbitrary analytic
function integrable in the interval $[0,1]$, if $\int^t_0S(t)\,dt$
is the sum function of the power series $\sum^\infty_{n=1}a_nt^n$,
let $f(x)\in{L^2}[a,b]$ be the sum function of the corresponding
cosine series, and $g(x)\in{L^2}[a,b]$ be that of the corresponding
sine series, namely
\begin{eqnarray*}
  \int^t_0\!\!S(t)\,dt &=& \sum^\infty_{n=1}a_nt^n,\qquad{t}\in\mathbb{R}^1,\;\;0\leq{t}\leq{r},\;\;0<r<+\infty, \\
  f(x) &=& \sum^\infty_{n=1}a_n\cos\frac{n\pi{x}}{c}, \\
  g(x) &=&
  \sum^\infty_{n=1}a_n\sin\frac{n\pi{x}}{c},\qquad{x}\in\mathbb{R}^1,\;\;a<x<b.
\end{eqnarray*}
Then we have the following mapping relations:

\parbox{11cm}{\begin{eqnarray*}\label{15}
f(x) &=&
\int^1_0\!\!S(\xi)\,d\xi-\frac{\pi}{c}\int^x_0\!\!\left.\sin\left(\frac{\pi{x}}{c}\frac{\partial}{\partial{z}}\right)
[S(e^z)\,e^z]\right|_{z=0}dx, \\
g(x)&=&\frac{\pi}{c}\int^x_0\!\!\left.\cos\left(\frac{\pi{x}}{c}\frac{\partial}{\partial{z}}\right)[S(e^z)\,e^z]\right|_{z=0}dx.
\end{eqnarray*}}\hfill\parbox{1cm}{\begin{eqnarray}\end{eqnarray}}
And the endpoints $a$ and $b$ of the interval $a<x<b$ are
non-analytical points (singularities) of Fourier series, which can
be uniquely determined by the detailed computation of the right-hand
side of (\ref{15}).

\textbf{Proof.} Combining Theorem 2 with Theorem 3 will lead us to
the proof.

Apparently, such summation method of Fourier series can be extended
into other trigonometric series.

\textbf{Theorem 5.}  Let $S(t)\in{J}$ be the sum function of power
series $\sum^\infty_{n=0}a_nt^n$, $f(x)$ be the sum function of the
corresponding $\sum^\infty_{n=0}a_n\cos(nx)\cos^nx$, $g(x)$ be the
sum function of the corresponding
$\sum^\infty_{n=0}a_n\sin(nx)\cos^nx$, namely
\begin{eqnarray*}
  S(t) &=& \sum^\infty_{n=0}a_nt^n,\qquad{t}\in\mathbb{R}^1,\;\;0\leq{t}\leq{r},\;\;0<r<+\infty, \\
  f(x) &=& \sum^\infty_{n=0}a_n\cos(nx)\cos^nx, \\
  g(x) &=&
  \sum^\infty_{n=0}a_n\sin(nx)\cos^nx,\qquad{x}\in\mathbb{R}^1,\;\;a<x<b,
\end{eqnarray*}
then we have the following mapping relations:

\parbox{11cm}{\begin{eqnarray*}\label{16}
 f(x) &=&
 \left.\cos\left(x\rho\frac{\partial}{\partial\rho}\right)S(\rho)\right|_{\rho=\cos{x}}=\;\cos\left(Y\frac{\partial}{\partial{X}}\right)S(X), \\
 g(x) &=&
 \left.\sin\left(x\rho\frac{\partial}{\partial\rho}\right)S(\rho)\right|_{\rho=\cos{x}}=\;\sin\left(Y\frac{\partial}{\partial{X}}\right)S(X).
 \end{eqnarray*}}\hfill\parbox{1cm}{\begin{eqnarray}\end{eqnarray}}
Where $X=\rho\cos{x}=\cos^2x$, $Y=\rho\sin{x}=\cos{x}\sin{x}$. And
the endpoints $a$ and $b$ of the interval $a<x<b$ are non-analytical
points (singularities) of trigonometric series, which can be
uniquely determined by the detailed computation of the right-hand
side of (\ref{16}).

\textbf{Proof.}  According to the Definition 5 in reference
\cite{bi10}, the following two expressions are obvious:
\begin{eqnarray*}
 f(x) &=&
 \left.\cos\left(x\rho\frac{\partial}{\partial\rho}\right)S(\rho)\right|_{\rho=\cos{x}}, \\
 g(x) &=&
 \left.\sin\left(x\rho\frac{\partial}{\partial\rho}\right)S(\rho)\right|_{\rho=\cos{x}},
 \end{eqnarray*}
 where $\cos\left(x\rho\frac{\partial}{\partial\rho}\right)$ and $\sin\left(x\rho\frac{\partial}{\partial\rho}\right)$
 are the abstract operators taking $\rho\frac{\partial}{\partial\rho}$ as the operator element. By using the Theorem 3 in reference \cite{bi10}, namely

 Let $\rho\in{R}^n$, $\theta\in{R}_n$, $X=(\rho_1\cos\theta_1,\ldots,\rho_n\cos\theta_n)$,
 $Y=(\rho_1\sin\theta_1,\ldots,\rho_n\sin\theta_n)$,
 then for an arbitrary analytic function $f(\rho)$, we have

\parbox{11cm}{\begin{eqnarray*}\label{17}
 \cos\left(\theta\rho\frac{\partial}{\partial\rho}\right)f(\rho) &=& \cos(Y\partial_X)f(X), \\
 \sin\left(\theta\rho\frac{\partial}{\partial\rho}\right)f(\rho) &=&
 \sin(Y\partial_X)f(X),
 \end{eqnarray*}}\hfill\parbox{1cm}{\begin{eqnarray}\end{eqnarray}}
where \[\theta\rho\frac{\partial}{\partial\rho}=
\left(\theta_1\rho_1\frac{\partial}{\partial\rho_1}+\cdots+\theta_n\rho_n\frac{\partial}{\partial\rho_n}\right),\]
we can obtain (\ref{16}) immediately.

\textbf{Corollary 1.} If the power series expansion of an analytic
function is unique, then according to the mapping relations between
analytic functions and periodic functions, the Fourier series
expansion of a periodic function is unique as well.

It is easy to prove the uniqueness of power series expansions of
analytic functions, thus Corollary 1 can actually derive the
uniqueness of Fourier series expansions of periodic functions
easier.

Theorem 2 can be directly extended into the multiple Fourier series
while the form of expressions remains essentially constant, namely

\textbf{Theorem 6.}  Let $S(t)\in{J}(\mathfrak{D})$ be the sum
function of power series $\sum^\infty_{n=0}a_nt^n$,
$f(x)\in{L^2}(\Omega)$ be the sum function of the corresponding
cosine series, and $g(x)\in{L^2}(\Omega)$ be that of the
corresponding sine series, namely
\begin{eqnarray*}
  S(t) &=& \sum^\infty_{n=0}a_nt^n,\qquad{t}\in\mathfrak{D}\subset\mathbb{R}^m. \\
  f(x) &=& \sum^\infty_{n=0}a_n\cos\frac{n\pi{x}}{c}, \\
  g(x) &=&
  \sum^\infty_{n=0}a_n\sin\frac{n\pi{x}}{c},\qquad{n}\in\mathbb{N}^m,\;\;x\in\Omega\subset\mathbb{R}^m,
\end{eqnarray*}
then we have the following mapping relations:

\parbox{11cm}{\begin{eqnarray*}\label{18}
                f(x) &=& \left.\cos\left(\frac{\pi{x}}{c}\frac{\partial}{\partial{z}}\right)S(e^z)\right|_{z=0}, \\
                g(x) &=&
                \left.\sin\left(\frac{\pi{x}}{c}\frac{\partial}{\partial{z}}\right)S(e^z)\right|_{z=0},
                \qquad(z\in\mathbb{R}^m).
              \end{eqnarray*}}\hfill\parbox{1cm}{\begin{eqnarray}\end{eqnarray}}
And the non-analytical points (singularities) on the border
$\partial\Omega$ can be uniquely determined by the detailed
computation of the right-hand side of (\ref{18}).

Here the signs appearing in these formulas should be interpreted as
the following universal abbreviation, namely
\[nx=n_1x_1+n_2x_2+\cdots+n_mx_m,\quad{a}_n=a_{n_1,n_2,\ldots,n_m},\quad{e}^z=\left(e^{z_1},e^{z_2},\ldots,e^{z_m}\right),\]
\[x\frac{\partial}{\partial{z}}=x_1\frac{\partial}{\partial{z_1}}+x_2\frac{\partial}{\partial{z_2}}+\cdots+
x_m\frac{\partial}{\partial{z_m}},\quad\sum^\infty_{n=0}=\sum^\infty_{n_1=0}\sum^\infty_{n_2=0}\cdots\sum^\infty_{n_m=0}.\]
And the meanings of other $t,x,t^n$ are the same with the universal
signs.

\textbf{Example 1.}  Trigonometric series: the sum function of
$\sum^\infty_{n=1}(1/n)\sin(nx)\cos^nx$ is $g(x)=\pi/2-x$, and its
non-analytical points are $x=0$ and $x=\pi$, namely
\begin{equation}\label{21}
\frac{\pi}{2}-x=\sum^\infty_{n=1}\frac{1}{n}\sin(nx)\cos^nx,\quad0<x<\pi.
\end{equation}

\textbf{Proof.}  By using Theorem 5, and the algorithms and basic
formulas in this paper, we have
\begin{eqnarray*}
  S(t) &=& \sum^\infty_{n=1}\frac{1}{n}t^n=-\ln(1-t),\quad{t}\in\mathbb{R}^1,\;|t|<1, \\
  g(t) &=&
  \sum^\infty_{n=1}\frac{1}{n}\sin(nx)\cos^nx,\quad{x}\in\mathbb{R}^1,\;a<x<b.
\end{eqnarray*}
\begin{eqnarray*}
  g(x) &=& \left.\sin\left(x\rho\frac{\partial}{\partial\rho}\right)S(\rho)\right|_{\rho=\cos{x}}
  =\sin\left(Y\frac{\partial}{\partial{X}}\right)S(X)=-\sin\left(Y\frac{\partial}{\partial{X}}\right)\ln(1-X) \\
       &=& \textrm{arccot}\frac{1-X}{Y}=
       \textrm{arccot}\frac{1-\cos^2x}{\sin{x}\cos{x}}=
       \textrm{arccot}\frac{\sin^2x}{\sin{x}\cos{x}}.
\end{eqnarray*}

When $\sin{x}\neq0$, namely $x\neq0$ and $x\neq\pi$ (thus $a=0$,
$b=\pi$), we have
\[g(x)=\textrm{arccot}\frac{\sin{x}}{\cos{x}}=\textrm{arccot}\tan{x}=\textrm{arccot}\cot(\frac{\pi}{2}-x)=\frac{\pi}{2}-x.\]

\textbf{Example 2.}  Fourier series: the sum function of
$\sum^\infty_{n=1}[(-1)^{n-1}/((3n-1)(3n+1))]\cos(3n\omega{t})$ is
$f(t)=(\sqrt{3}\pi/9)\cos\omega{t}-1/2$, and its non-analytical
point is $|\omega{t}|=\pi/3$, namely
\begin{equation}\label{22}
    \frac{\sqrt{3}}{9}\pi\cos\omega{t}-\frac{1}{2}=\sum^\infty_{n=1}\frac{(-1)^{n-1}}{(3n-1)(3n+1)}\cos(3n\omega{t}),\quad
    -\frac{\pi}{3}<\omega{t}<\frac{\pi}{3}.
\end{equation}

\textbf{Proof.}  By using Theorem 2, and the algorithms and basic
formulas in this paper, let
\begin{eqnarray*}
  S(x) &=& \sum^\infty_{n=1}(-1)^{n-1}\frac{x^{3n}}{(3n-1)(3n+1)},\quad{x}\in\mathbb{R}^1,\;|x|<1, \\
  f(t) &=&
  \sum^\infty_{n=1}\frac{(-1)^{n-1}}{(3n-1)(3n+1)}\cos(3n\omega{t}),\quad{t}\in\mathbb{R}^1,\;\,a<\omega{t}<b.
\end{eqnarray*}

As it is difficult to obtain the sum function $S(x)$ directly, we
can use the operators $\frac{d}{dx}(x\cdot)$ and
$\frac{d}{dx}(\frac{1}{x}\cdot)$ to transform the power series into
the series familiar to us, then we can obtain $S(x)$, and then
$f(t)$, namely
\[\frac{d}{dx}\left(\frac{1}{x}\frac{d}{dx}(xS(x))\right)=\sum^\infty_{n=1}(-1)^{n-1}x^{3n-2}=\frac{x}{1+x^3}.\]

To describe concisely, the following results will be given directly:
\begin{eqnarray*}
  S(x) &=& \frac{1}{x}\int^x_0\!\!x\,dx\int^x_0\!\!\frac{x}{1+x^3}\,dx=\left(\frac{x}{12}-\frac{1}{12x}\right)\ln\left(x^2-x+1\right)\\
       & &
       -\left(\frac{x}{6}-\frac{1}{6x}\right)\ln(1+x)+
       \left(\frac{x}{4}+\frac{1}{4x}\right)\frac{2}{\sqrt{3}}\left(\arctan\frac{2x-1}{\sqrt{3}}+\frac{\pi}{6}\right)-\frac{1}{2}.
\end{eqnarray*}
\begin{eqnarray*}
f(t) &=& \left.\cos\left(\omega{t}\frac{\partial}{\partial{z}}\right)S(e^z)\right|_{z=0} \\
&=&
\frac{1}{6}\left.\cos\left(\omega{t}\frac{\partial}{\partial{z}}\right)\left[\,\sinh{z}\ln\left(e^{2z}-e^z+1\right)\right]\right|_{z=0}\\
& &-\frac{1}{3}\left.\cos\left(\omega{t}\frac{\partial}{\partial{z}}\right)[\,\sinh{z}\ln\left(1+e^z\right)]\right|_{z=0} \\
& &
 +\frac{1}{\sqrt{3}}\left.\cos\left(\omega{t}\frac{\partial}{\partial{z}}\right)\left[\,\cosh{z}\left(\arctan\frac{2e^z-1}{\sqrt{3}}
  +\frac{\pi}{6}\right)\right]\right|_{z=0}-\frac{1}{2}\\
&=&
-\frac{1}{6}\sin\omega{t}\,\textrm{arccot}\frac{(2\cos\omega{t}-1)\cos\omega{t}}{(2\cos\omega{t}-1)\sin\omega{t}}
+\frac{1}{3}\sin\omega{t}\,\textrm{arccot}\frac{\cos^2(\omega{t}/2)}{\sin(\omega{t}/2)\cos(\omega{t}/2)} \\
& &
+\frac{1}{2\sqrt{3}}\cos\omega{t}\arctan\frac{\sqrt{3}(2\cos\omega{t}-1)}{2\cos\omega{t}-1}+\frac{\pi}{6\sqrt{3}}\cos\omega{t}
-\frac{1}{2}.
\end{eqnarray*}

When $2\cos\omega{t}-1\neq0$, and $\cos(\omega{t}/2)\neq0$, namely
$|\omega{t}|\neq\pi/3$ and $|\omega{t}|\neq\pi$ , we have
\begin{eqnarray*}
  f(t) &=& -\frac{1}{6}\sin\omega{t}\,\textrm{arccot}\cot\omega{t}+\frac{1}{3}\sin\omega{t}\,\textrm{arccot}\cot\frac{\omega{t}}{2} \\
       & &
       +\frac{1}{2\sqrt{3}}\cos\omega{t}\arctan\sqrt{3}+\frac{\pi}{6\sqrt{3}}\cos\omega{t}-\frac{1}{2}
       =\frac{\sqrt{3}}{9}\pi\cos\omega{t}-\frac{1}{2}.
\end{eqnarray*}
There are four non-analytical points in the interval $[-\pi,\pi]$:
$\omega{t}=-\pi,-\pi/3,\pi/3,\pi$, thus $a=-\pi/3,\;b=\pi/3$.
Therefore, Example 2 is proved.

\section{The Zeta function of odd variables}
\noindent

\textbf{Definition 1.}  Let $S_0(t)$ be a function analytic in the
neighborhood of $t=0$ and
\[S_0(t)=\sum^\infty_{n=1}a_nt^n,\qquad|t|<r,\quad0<r<+\infty,\]
then $S_m(t)$ is defined as
\begin{equation}\label{23}
S_m(t)=\underbrace{\int^t_0\frac{dt}{t}\cdots}_m\int^t_0S_0(t)\,\frac{dt}{t}=\sum^\infty_{n=1}a_n\frac{t^n}{n^m}.
\end{equation}
Apparently $S_m(1)$ is the sum function of the Dirichlet series
taking $m$ as the variable.

\textbf{Lemma 1.}  According to (\ref{23}), $S_m(t)$ satisfies the
following recurrence relation:
\begin{equation}\label{24}
\int^t_0S_{m-1}(t)\,\frac{dt}{t}=S_m(t).
\end{equation}

\textbf{Lemma 2.}  The sum function $S_m(1)$ has the following
recurrence property:
\begin{equation}\label{25}
\left.\cos\left(\frac{\pi{x}}{c}\frac{\partial}{\partial{z}}\right)S_m(e^z)\right|_{z=0}=S_m(1)
-\frac{\pi}{c}\int^x_0\left.\sin\left(\frac{\pi{x}}{c}\frac{\partial}{\partial{z}}\right)S_{m-1}(e^z)\right|_{z=0}dx.
\end{equation}

\textbf{Proof.}  Taking $S(x)=S_{m-1}(x)/x$ in (\ref{13}), it is
proved by using Lemma 1.

Similarly,

\textbf{Lemma 3.}  $S_m(t)$ has the following recurrence property:
\begin{equation}\label{26}
\left.\sin\left(\frac{\pi{x}}{c}\frac{\partial}{\partial{z}}\right)S_m(e^z)\right|_{z=0}=
\frac{\pi}{c}\int^x_0\left.\cos\left(\frac{\pi{x}}{c}\frac{\partial}{\partial{z}}\right)S_{m-1}(e^z)\right|_{z=0}dx.
\end{equation}

\textbf{Theorem 7.}  The sum function $S_m(1)$ has the following
property:

\parbox{11cm}{\begin{eqnarray*}\label{27}
                & &
                \left.\cos\left(\frac{\pi{x}}{c}\frac{\partial}{\partial{z}}\right)S_{m-1}(e^z)\right|_{z=0}\\
                &=& \sum^{r-1}_{k=0}(-1)^k\frac{1}{(2k)!}\left(\frac{\pi{x}}{c}\right)^{2k}S_{m-1-2k}(1) \\
                & &
                +\,(-1)^r\left(\frac{\pi}{c}\right)^{2r-1}
\underbrace{\int^x_0dx\cdots}_{2r-1}\int^x_0\left.\sin\left(\frac{\pi{x}}{c}\frac{\partial}{\partial{z}}\right)S_{m-2r}(e^z)\right|_{z=0}dx.
              \end{eqnarray*}}\hfill\parbox{1cm}{\begin{eqnarray}\end{eqnarray}}

\textbf{Proof.}  We can use the mathematical induction to prove it.
According to Theorem 3, it is obviously tenable when $r=1$ in
(\ref{27}). Now we inductively hypothesize that it is tenable when
$r=K$, namely
\begin{eqnarray*}
 & & \left.\cos\left(\frac{\pi{x}}{c}\frac{\partial}{\partial{z}}\right)S_{m-1}(e^z)\right|_{z=0}\\
 &=& \sum^{K-1}_{k=0}(-1)^k\frac{1}{(2k)!}\left(\frac{\pi{x}}{c}\right)^{2k}S_{m-1-2k}(1) \\
 & & +\,(-1)^K\left(\frac{\pi}{c}\right)^{2K-1}
\underbrace{\int^x_0dx\cdots}_{2K-1}\int^x_0\left.\sin\left(\frac{\pi{x}}{c}\frac{\partial}{\partial{z}}\right)S_{m-2K}(e^z)\right|_{z=0}dx.
\end{eqnarray*}
Using Lemma 3 and 2 respectively, then the above expression turns
into
\begin{eqnarray*}
 & & \left.\cos\left(\frac{\pi{x}}{c}\frac{\partial}{\partial{z}}\right)S_{m-1}(e^z)\right|_{z=0}\\
 &=& \sum^{K-1}_{k=0}(-1)^k\frac{1}{(2k)!}\left(\frac{\pi{x}}{c}\right)^{2k}S_{m-1-2k}(1) \\
 & & +\,(-1)^K\left(\frac{\pi}{c}\right)^{2K}
\underbrace{\int^x_0dx\cdots}_{2K}\int^x_0\left.\cos\left(\frac{\pi{x}}{c}\frac{\partial}{\partial{z}}\right)S_{m-2K-1}(e^z)\right|_{z=0}dx\\
 &=&
 \sum^{K-1}_{k=0}(-1)^k\frac{1}{(2k)!}\left(\frac{\pi{x}}{c}\right)^{2k}S_{m-1-2k}(1)+
 (-1)^K\left(\frac{\pi}{c}\right)^{2K}S_{m-2K-1}(1)\frac{x^{2K}}{(2K)!}\\
 & & +\,(-1)^{K+1}\left(\frac{\pi}{c}\right)^{2K+1}
\underbrace{\int^x_0dx\cdots}_{2K+1}\int^x_0\left.\sin\left(\frac{\pi{x}}{c}\frac{\partial}{\partial{z}}\right)S_{m-2K-2}(e^z)\right|_{z=0}dx\\
 &=&
 \sum^K_{k=0}(-1)^k\frac{1}{(2k)!}\left(\frac{\pi{x}}{c}\right)^{2k}S_{m-1-2k}(1)\\
 & & +\,(-1)^{K+1}\left(\frac{\pi}{c}\right)^{2K+1}
\underbrace{\int^x_0dx\cdots}_{2K+1}\int^x_0\left.\sin\left(\frac{\pi{x}}{c}\frac{\partial}{\partial{z}}\right)S_{m-2K-2}(e^z)\right|_{z=0}dx.
\end{eqnarray*}
Thus it is tenable when $r=K+1$, and then Theorem 7 is proved.

Similarly we can prove the following theorems:

\textbf{Theorem 8.}  The sum function $S_m(1)$ has the following
property:

\parbox{11cm}{\begin{eqnarray*}\label{28}
                & &
                \left.\sin\left(\frac{\pi{x}}{c}\frac{\partial}{\partial{z}}\right)S_{m-1}(e^z)\right|_{z=0}\\
                &=& \sum^{r-1}_{k=1}(-1)^{k-1}\frac{1}{(2k-1)!}\left(\frac{\pi{x}}{c}\right)^{2k-1}S_{m-2k}(1) \\
                & &
                +\,(-1)^{r-1}\left(\frac{\pi}{c}\right)^{2r-1}
\underbrace{\int^x_0dx\cdots}_{2r-1}\int^x_0\left.\cos\left(\frac{\pi{x}}{c}\frac{\partial}{\partial{z}}\right)S_{m-2r}(e^z)\right|_{z=0}dx.
              \end{eqnarray*}}\hfill\parbox{1cm}{\begin{eqnarray}\end{eqnarray}}

\textbf{Theorem 9.}  The sum function $S_m(1)$ has the following
property:

\parbox{11cm}{\begin{eqnarray*}\label{29}
                & &
                \left.\sin\left(\frac{\pi{x}}{c}\frac{\partial}{\partial{z}}\right)S_{m-1}(e^z)\right|_{z=0}\\
                &=& \sum^r_{k=1}(-1)^{k-1}\frac{1}{(2k-1)!}\left(\frac{\pi{x}}{c}\right)^{2k-1}S_{m-2k}(1) \\
                & &
                +\,(-1)^r\left(\frac{\pi}{c}\right)^{2r}
\underbrace{\int^x_0dx\cdots}_{2r}\int^x_0\left.\sin\left(\frac{\pi{x}}{c}\frac{\partial}{\partial{z}}\right)S_{m-2r-1}(e^z)\right|_{z=0}dx.
              \end{eqnarray*}}\hfill\parbox{1cm}{\begin{eqnarray}\end{eqnarray}}

\textbf{Theorem 10.}  The sum function $S_m(1)$ has the following
property:

\parbox{11cm}{\begin{eqnarray*}\label{30}
                & &
                \left.\cos\left(\frac{\pi{x}}{c}\frac{\partial}{\partial{z}}\right)S_{m-1}(e^z)\right|_{z=0}\\
                &=& \sum^{r-1}_{k=0}(-1)^k\frac{1}{(2k)!}\left(\frac{\pi{x}}{c}\right)^{2k}S_{m-2k-1}(1) \\
                & &
                +\,(-1)^r\left(\frac{\pi}{c}\right)^{2r}
\underbrace{\int^x_0dx\cdots}_{2r}\int^x_0\left.\cos\left(\frac{\pi{x}}{c}\frac{\partial}{\partial{z}}\right)S_{m-2r-1}(e^z)\right|_{z=0}dx.
              \end{eqnarray*}}\hfill\parbox{1cm}{\begin{eqnarray}\end{eqnarray}}

\textbf{Lemma 4.} In the interval $(0,2c)$, we have:

\parbox{11cm}{\begin{eqnarray*}\label{31}
 \left.\sin\left(\frac{\pi{x}}{c}\frac{\partial}{\partial{z}}\right)(-\ln(1-e^z))\right|_{z=0}&=&
 \frac{\pi}{2}-\frac{\pi{x}}{2c},\quad0<x<2c.\\
 \left.\sin\left(\frac{\pi{x}}{c}\frac{\partial}{\partial{z}}\right)\ln(1+e^z)\right|_{z=0}&=&
 \frac{\pi{x}}{2c},\quad|x|<c.\\
 \left.\cos\left(\frac{\pi{x}}{c}\frac{\partial}{\partial{z}}\right)\arctan{e^z}\right|_{z=0}&=& \frac{\pi}{4},\quad|x|<c/2.
\end{eqnarray*}}\hfill\parbox{1cm}{\begin{eqnarray}\end{eqnarray}}

\textbf{Proof.} By using the algorithms and basic formulas in this
paper, we have
\begin{eqnarray*}
   & &\left.\sin\left(\frac{\pi{x}}{c}\frac{\partial}{\partial{z}}\right)(-\ln(1-e^z))\right|_{z=0}\,=\,
   -\left.\sin\left(Y\frac{\partial}{\partial{X}}\right)\ln{X}\right|_{z=0}=-\left.\textrm{arccot}\frac{X}{Y}\right|_{z=0}\\
   &=&
   -\textrm{arccot}\frac{1-\cos(\pi{x}/c)}{-\sin(\pi{x}/c)}=\textrm{arccot}\frac{\sin^2(\pi x/(2c))}{\sin(\pi x/(2c))\cos(\pi
   x/(2c))}.
\end{eqnarray*}
When $\sin(\pi{x}/(2c))\neq0$ or $x\neq0,\;x\neq2c$, the above
expression turns into
\[\left.\sin\left(\frac{\pi{x}}{c}\frac{\partial}{\partial{z}}\right)(-\ln(1-e^z))\right|_{z=0}=
\textrm{arccot}\tan\frac{\pi{x}}{2c}=\frac{\pi}{2}-\frac{\pi{x}}{2c}.\]
Similarly we have
\begin{eqnarray*}
   & &\left.\sin\left(\frac{\pi{x}}{c}\frac{\partial}{\partial{z}}\right)\ln(1+e^z)\right|_{z=0}\,=\,
   \left.\sin\left(Y\frac{\partial}{\partial{X}}\right)\ln{X}\right|_{z=0}=\left.\textrm{arccot}\frac{X}{Y}\right|_{z=0}\\
   &=&
   \textrm{arccot}\frac{1+\cos(\pi{x}/c)}{\sin(\pi{x}/c)}=\textrm{arccot}\frac{\cos^2(\pi x/(2c))}{\sin(\pi x/(2c))\cos(\pi
   x/(2c))}.
\end{eqnarray*}
When $\cos(\pi{x}/(2c))\neq0$ or $|x|\neq{c}$, the above expression
turns into
\[\left.\sin\left(\frac{\pi{x}}{c}\frac{\partial}{\partial{z}}\right)\ln(1+e^z)\right|_{z=0}=
\textrm{arccot}\cot\frac{\pi{x}}{2c}=\frac{\pi{x}}{2c}.\]
\begin{eqnarray*}
   & & \left.\cos\left(\frac{\pi{x}}{c}\frac{\partial}{\partial{z}}\right)\arctan{e^z}\right|_{z=0}\,=\,
   \left.\cos\left(Y\frac{\partial}{\partial{X}}\right)\arctan{X}\right|_{z=0}\\
   &=&\left.\frac{1}{2}\arctan\frac{2X}{1-(X^2+Y^2)}\right|_{z=0}\,=\,
   \frac{1}{2}\arctan\frac{2\cos(\pi{x}/c)}{1-\left(\cos^2(\pi{x}/c)+\sin^2(\pi{x}/c)\right)}.
\end{eqnarray*}
When $\cos(\pi{x}/c)\neq0$ or $|x|\neq{c/2}$, the above expression
turns into
\[\left.\cos\left(\frac{\pi{x}}{c}\frac{\partial}{\partial{z}}\right)\arctan{e^z}\right|_{z=0}=\frac{1}{2}\arctan\infty=\frac{\pi}{4}.\]

\textbf{Theorem 11.}  In the interval $(0,2c)$, we have the
following Fourier series expressions:

\parbox{11cm}{\begin{eqnarray*}\label{32}
 & &
 \sum^{\infty}_{n=1}\frac{1}{n^{2r}}\cos\frac{n\pi{x}}{c}\\
 &=& \sum^{r-1}_{k=0}(-1)^{k}\frac{1}{(2k)!}\left(\frac{\pi{x}}{c}\right)^{2k}\zeta(2r-2k)\\
 & & +\,(-1)^{r}\left(\frac{\pi}{2}\frac{1}{(2r-1)!}\left(\frac{\pi{x}}{c}\right)^{2r-1}
-\frac{1}{2}\frac{1}{(2r)!}\left(\frac{\pi{x}}{c}\right)^{2r}\right).
\end{eqnarray*}}\hfill\parbox{1cm}{\begin{eqnarray}\end{eqnarray}}
\parbox{11cm}{\begin{eqnarray*}\label{33}
 & &
 \sum^{\infty}_{n=1}\frac{1}{n^{2r+1}}\sin\frac{n\pi{x}}{c}\\
 &=& \sum^{r}_{k=1}(-1)^{k-1}\frac{1}{(2k-1)!}\left(\frac{\pi{x}}{c}\right)^{2k-1}\zeta(2r+2-2k)\\
 & &
+\,(-1)^{r}\left(\frac{\pi}{2}\frac{1}{(2r)!}\left(\frac{\pi{x}}{c}\right)^{2r}
-\frac{1}{2}\frac{1}{(2r+1)!}\left(\frac{\pi{x}}{c}\right)^{2r+1}\right).
\end{eqnarray*}}\hfill\parbox{1cm}{\begin{eqnarray}\end{eqnarray}}

\textbf{Proof.}  According to Definition 1, if $S_0(t)=-\ln(1-t)$,
then
\[S_{m-1}(t)=\sum^\infty_{n=1}\frac{t^n}{n^m},\quad{S_{m-1}(1)}=\zeta(m).\]
In Theorem 7, let $m=2r$ , and by using Theorem 2, we have
\begin{eqnarray*}
                & &
                 \sum^{\infty}_{n=1}\frac{1}{n^{2r}}\cos\frac{n\pi{x}}{c}\\
                &=& \sum^{r-1}_{k=0}(-1)^k\frac{1}{(2k)!}\left(\frac{\pi{x}}{c}\right)^{2k}S_{2r-1-2k}(1) \\
                & &
                +\,(-1)^r\left(\frac{\pi}{c}\right)^{2r-1}
\underbrace{\int^x_0dx\cdots}_{2r-1}\int^x_0\left.\sin\left(\frac{\pi{x}}{c}\frac{\partial}{\partial{z}}\right)(-\ln(1-e^z))\right|_{z=0}dx.
              \end{eqnarray*}
Where $S_{2r-1-2k}(1)=\zeta(2r-2k)$, by using Lemma 4 we have
(\ref{32}).

Let $m=2r+1$ in Theorem 9, using Theorem 2, and considering
$S_{2r+1-2k}(1)=\zeta(2r+2-2k)$ and Lemma 4, we have (\ref{33}).

\textbf{Corollary 2.}  In Theorem 11, when $r\geq1$, the Fourier
series expression is tenable at the endpoints of the interval.

\textbf{Proof.}  Observing the summation formulas of trigonometric
series given by Theorem 11, if we use the integral operator
$\int^x_0dx$ to integrate both sides of (\ref{32}), we can obtain
(\ref{33}). It means that, for any point $x$ in the interval of
convergence of summation formulas of trigonometric series, if $x=0$
is within the interval or at the endpoints of it, then the termwise
integration of the trigonometric series by the integral operator
$\int^x_0dx$ converges uniformly to the integration of the sum
function. Actually, in the theory of Fourier series, according to
the integrability, if there is a certain summation formula of
trigonometric series in the open interval $a<x<b$, then the termwise
integration of the trigonometric series by the integral operator
$\int^x_adx,\;x\in[a,b]$ converges uniformly to the integration of
the sum function. In other words, for a summation formula of
trigonometric series in $a<x<b$, if we use the integral operator
$\int^x_adx,\;x\in[a,b]$ to integrate both sides of the equation,
then the new summation formula of trigonometric series obtained is
definitely tenable in the closed interval $a\leq{x}\leq{b}$. Thus
Inference 3 is proved.

\textbf{Theorem 12.}  The Riemann Zeta function $\zeta(2n)$
satisfies the following recurrence formula:
\begin{equation}\label{34}
    \sum^{r-1}_{k=0}(-1)^k\frac{\pi^{2k}}{(2k+1)!}\zeta(2r-2k)=(-1)^{r-1}\frac{\pi^{2r}}{(2r+1)!}r.
\end{equation}

\textbf{Proof}  Let $x=c$ in the second expression of Theorem 11,
then it can be proved.

\textbf{Theorem 13.}  The Riemann Zeta function $\zeta(2n)$
satisfies the following recursion formula:
\begin{equation}\label{35}
    \sum^{r-1}_{k=0}(-1)^k\frac{(2\pi)^{2k}}{(2k+1)!}\zeta(2r-2k)=(-1)^{r-1}\frac{(2\pi)^{2r}}{4(2r+1)!}(2r-1).
\end{equation}

\textbf{Proof.}  Let $x=2c$ in the second expression of Theorem 11,
then it can be proved.

\textbf{Theorem 14.}  In the interval $(0,2c)$, we have the
following Fourier series expressions related to $\zeta(2n+1)$:

\parbox{11cm}{\begin{eqnarray*}\label{36}
 & &
 \sum^\infty_{n=1}\frac{1}{n^{2r+1}}\cos\frac{n\pi{x}}{c}\\
 &=&
 \sum^{r-1}_{k=0}(-1)^{k}\frac{1}{(2k)!}\left(\frac{\pi{x}}{c}\right)^{2k}\zeta(2r+1-2k)\\
 & & +\,(-1)^{r-1}\left(\frac{\pi}{c}\right)^{2r}\underbrace{\int^x_0dx\cdots}_{2r}\int^x_0\ln\left(2\sin\frac{\pi{x}}{2c}\right)\,dx.
\end{eqnarray*}}\hfill\parbox{1cm}{\begin{eqnarray}\end{eqnarray}}

\textbf{Proof.} Taking $m=2r+1,\;S_0(t)=-\ln(1-t)$ in Theorem 10,
considering it is in the interval $(0,2c)$, and by using the
algorithms and basic formulas in this paper, we have
\begin{eqnarray*}
   & & \left.\cos\left(\frac{\pi{x}}{c}\frac{\partial}{\partial{z}}\right)(-\ln(1-e^z))\right|_{z=0} \\
   &=&
   -\left.\cos\left(Y\frac{\partial}{\partial{X}}\right)\ln{X}\right|_{z=0}=-\left.\frac{1}{2}\ln(X^2+Y^2)\right|_{z=0}\\
   &=&
   -\frac{1}{2}\ln\left(\left(1-\cos\frac{\pi{x}}{c}\right)^2+\sin^2\frac{\pi{x}}{c}\right)=-\ln\left(2\sin\frac{\pi{x}}{2c}\right),
\end{eqnarray*}
and $S_{2r-2k}(1)=\zeta(2r+1-2k)$, thus it is proved by using
Theorem 2.

\textbf{Lemma 5.}  Let $m$ be an arbitrary natural number, then for
$\ln{x}$ we have the following integral formula:
\begin{equation}\label{37}
\underbrace{\int^x_0dx\cdots}_m\int^x_0\ln{x}\,dx=\frac{x^m}{m!}\left(\ln{x}-H_m\right).
\end{equation}
Where $H_m$ are the harmonic numbers.

\textbf{Theorem 15.}  Let $r\in\mathbb{N}$ be an arbitrary natural
number, and $B_k^*$ be the Bernoulli numbers, $B_k^*$ satisfies the
following recurrence formula:
\begin{equation}\label{01}
\sum^{r-1}_{k=0}(-1)^k{2r+1\choose{2k+1}}B_{k+1}^*=\frac{1}{2},
\end{equation}
then the Riemann Zeta function $\zeta(2n+1)$ can be recursively
determined by the following recurrence formula, namely

\parbox{11cm}{\begin{eqnarray*}\label{38}
 \zeta(2r+1)&=& \frac{2^{4r+1}}{2^{4r+1}+2^{2r}-1}\sum^{r-1}_{k=1}(-1)^{k-1}\frac{1}{(2k)!}\left(\frac{\pi}{2}\right)^{2k}\zeta(2r+1-2k)\\
 & &
 +\,(-1)^{r-1}\frac{2^{2r+1}\pi^{2r}}{(2^{4r+1}+2^{2r}-1)(2r)!}\left(H_{2r}-\ln\frac{\pi}{2}\right)\\
 & &
 +\,(-1)^{r-1}\frac{2^{2r}\pi^{2r}}{2^{4r+1}+2^{2r}-1}\sum^\infty_{k=1}\frac{1}{k}\frac{B_k^*}{(2r+2k)!}\left(\frac{\pi}{2}\right)^{2k},
\end{eqnarray*}}\hfill\parbox{1cm}{\begin{eqnarray}\end{eqnarray}}
or, equivalently,

\parbox{11cm}{\begin{eqnarray*}\label{40}
 \zeta(2r+1)&=& \frac{2^{4r+1}}{2^{4r+1}+2^{2r}-1}\sum^{r-1}_{k=1}(-1)^{k-1}\frac{1}{(2k)!}\left(\frac{\pi}{2}\right)^{2k}\zeta(2r+1-2k)\\
 & &
 +\,(-1)^{r-1}\frac{2^{2r+1}\pi^{2r}}{(2^{4r+1}+2^{2r}-1)(2r)!}\left(H_{2r}-\ln\frac{\pi}{2}\right)\\
 & & +\,(-1)^{r-1}\frac{2^{2r+1}\pi^{2r}}{2^{4r+1}+2^{2r}-1}\sum^\infty_{k=1}\frac{(2k)!}{4^{2k}(2r+2k)!k}\zeta(2k).
\end{eqnarray*}}\hfill\parbox{1cm}{\begin{eqnarray}\end{eqnarray}}

\textbf{Proof.}  According to Theorem 14 and Lemma 5, in the
interval $(0,2c)$ (if $r\geq1$, then in the closed interval
$[0,2c]$) we have
\begin{eqnarray*}
   & & \sum^\infty_{n=1}\frac{1}{n^{2r+1}}\cos\frac{n\pi{x}}{c} \\
   &=& \sum^{r-1}_{k=0}(-1)^{k}\frac{1}{(2k)!}\left(\frac{\pi{x}}{c}\right)^{2k}\zeta(2r+1-2k)
   +\frac{(-1)^{r-1}}{(2r)!}\left(\frac{\pi{x}}{c}\right)^{2r}\ln2 \\
   & & +\,(-1)^{r-1}\left(\frac{\pi}{c}\right)^{2r}\underbrace{\int^x_0dx\cdots}_{2r}\int^x_0\left(\ln\frac{\pi{x}}{2c}
   +\ln\frac{\sin\frac{\pi{x}}{2c}}{\frac{\pi{x}}{2c}}\right)\,dx \\
   &=& \sum^{r-1}_{k=0}(-1)^{k}\frac{1}{(2k)!}\left(\frac{\pi{x}}{c}\right)^{2k}\zeta(2r+1-2k) \\
   & &
   +\,(-1)^{r}\frac{1}{(2r)!}\left(\frac{\pi{x}}{c}\right)^{2r}\left(H_{2r}-\ln\frac{\pi{x}}{c}\right)\\
   & &+\,(-1)^{r-1}\left(\frac{\pi}{c}\right)^{2r}\underbrace{\int^x_0dx\cdots}_{2r}\int^x_0\ln\frac{\sin\frac{\pi{x}}{2c}}{\frac{\pi{x}}{2c}}\,dx\\
   &=&\sum^{r-1}_{k=0}(-1)^{k}\frac{1}{(2k)!}\left(\frac{\pi{x}}{c}\right)^{2k}\zeta(2r+1-2k) \\
   & &
   +\,(-1)^{r}\frac{1}{(2r)!}\left(\frac{\pi{x}}{c}\right)^{2r}\left(H_{2r}-\ln\frac{\pi{x}}{c}\right)\\
   & &+\,(-1)^{r-1}\left(\frac{\pi}{c}\right)^{2r}\underbrace{\int^x_0dx\cdots}_{2r}\int^x_0
   -\sum^\infty_{k=1}\frac{1}{2k}\frac{2^{2k}B_k^*}{(2k)!}\left(\frac{\pi{x}}{2c}\right)^{2k}\,dx\\
   &=&
   \sum^{r-1}_{k=0}(-1)^{k}\frac{1}{(2k)!}\left(\frac{\pi{x}}{c}\right)^{2k}\zeta(2r+1-2k) \\
   & &+\,(-1)^{r}\frac{1}{(2r)!}\left(\frac{\pi{x}}{c}\right)^{2r}\left(H_{2r}-\ln\frac{\pi{x}}{c}\right)\\
   & &+\,(-1)^r\sum^\infty_{k=1}\frac{1}{2k}\frac{B_k^*}{(2r+2k)!}\left(\frac{\pi{x}}{c}\right)^{2r+2k}.
\end{eqnarray*}
Let $x=c/2$, considering
$\cos\frac{(2k-1)\pi}{2}=\sin{k\pi}=0,\;\cos\frac{2k\pi}{2}=\cos{k\pi}=(-1)^k$
we have
\begin{eqnarray*}
  -\frac{1}{2^{2r+1}}\eta(2r+1)&=& \sum^{r-1}_{k=0}(-1)^{k}\frac{1}{(2k)!}\left(\frac{\pi}{2}\right)^{2k}\zeta(2r+1-2k) \\
   & & +\,(-1)^{r}\frac{1}{(2r)!}\left(\frac{\pi}{2}\right)^{2r}\left(H_{2r}-\ln\frac{\pi}{2}\right) \\
   & & +\,(-1)^r\sum^\infty_{k=1}\frac{1}{2k}\frac{B_k^*}{(2r+2k)!}\left(\frac{\pi}{2}\right)^{2r+2k}.
\end{eqnarray*}
Where $\eta(2r+1)$ is the Dirichlet Eta function, defined usually by
the Dirichlet series:
\begin{equation}\label{02}
\eta(s)=\sum^\infty_{n=1}\frac{(-1)^{n-1}}{n^s},\qquad(\Re(s)>0).
\end{equation}
Apparently
\begin{equation}\label{39}
  \eta(2r+1)=\frac{2^{2r}-1}{2^{2r}}\zeta(2r+1),\quad{r}\in\mathbb{N}.
\end{equation}
Substituting (\ref{39}) into the above expression, then it is
proved.

\textbf{Example 3.}  In Theorem 15, taking $r=1$, we have
\begin{equation}\label{03}
\zeta(3)=\frac{6\pi^2}{35}-\frac{4\pi^2}{35}\ln\frac{\pi}{2}+\frac{\pi^2}{35}\sum^\infty_{k=1}\frac{B_k^*\pi^{2k}}{k4^{k-1}(2k+2)!}.
\end{equation}

In general, of which a useful generalization is

\textbf{Theorem 16.} We have the following relation between the
Riemann Zeta functions $\zeta(2r+1)$ and
$\zeta(2r),\;(r\in\mathbb{N})$:

\parbox{11cm}{\begin{eqnarray*}\label{36'}
 \sum^\infty_{n=1}\frac{1}{n^{2r+1}}\cos\frac{n\pi{x}}{c}
 &=& \sum^{r-1}_{k=0}(-1)^k\frac{1}{(2k)!}\left(\frac{\pi{x}}{c}\right)^{2k}\zeta(2r+1-2k)\\
 & &
 +\,(-1)^r\frac{1}{(2r)!}\left(\frac{\pi{x}}{c}\right)^{2r}\left(H_{2r}-\ln\frac{\pi{x}}{c}\right)\\
 & & +\,(-1)^r\left(\frac{\pi{x}}{c}\right)^{2r}\sum^\infty_{k=1}\frac{(2k-1)!\,\zeta(2k)}{2^{2k-1}(2r+2k)!}\left(\frac{x}{c}\right)^{2k},\quad0\leq{x}\leq2c.
\end{eqnarray*}}\hfill\parbox{1cm}{\begin{eqnarray}\end{eqnarray}}

In 1999, by using another method in reference \cite{sr}, H. M.
Srivastava has obtained the following results ($n\in\mathbb{N}$):

\parbox{11cm}{\begin{eqnarray*}\label{40'}
 \zeta(2n+1)&=& (-1)^{n-1}\frac{2(2\pi)^{2n}}{2^{4n+1}+2^{2n}-1}\left[\frac{H_{2n}-\ln(\frac{1}{2}\pi)}{(2n)!}\right.\\
 & &
 +\left.\sum^{n-1}_{k=1}\frac{(-1)^k}{(2n-2k)!}\frac{\zeta(2k+1)}{(\frac{1}{2}\pi)^{2k}}+2\sum^\infty_{k=1}\frac{(2k-1)!}{(2n+2k)!}\frac{\zeta(2k)}{4^{2k}}\right].
\end{eqnarray*}}\hfill\parbox{1cm}{\begin{eqnarray}\end{eqnarray}}

\parbox{11cm}{\begin{eqnarray*}\label{40''}
 \zeta(2n+1)&=& (-1)^{n-1}\frac{2(2\pi)^{2n}}{3^{2n}(2^{2n}+1)+2^{2n}-1}\left[\frac{H_{2n}-\ln(\frac{1}{3}\pi)}{(2n)!}\right.\\
 & &
 +\left.\sum^{n-1}_{k=1}\frac{(-1)^k}{(2n-2k)!}\frac{\zeta(2k+1)}{(\frac{1}{3}\pi)^{2k}}+2\sum^\infty_{k=1}\frac{(2k-1)!}{(2n+2k)!}\frac{\zeta(2k)}{6^{2k}}\right].
\end{eqnarray*}}\hfill\parbox{1cm}{\begin{eqnarray}\end{eqnarray}}
Obviously (\ref{40}) and (\ref{40'}) are identical. By using Theorem
16 and properties of the Hurwitz Zeta function, we also obtain the
(\ref{40''}).

The (Hurwitz's) generalized Zeta function $\zeta(s,a)$ and the Lerch
transcendent $\Phi(z,s,a)$ are usually defined by
\begin{equation}\label{yd1}
\zeta(s,a)=\sum^\infty_{n=0}\frac{1}{(n+a)^s},\quad(\Re(s)>1;\;a\neq0,-1,-2,\cdots)
\end{equation}
and
\begin{equation}\label{yd2}
\Phi(z,s,a)=\sum^\infty_{n=0}\frac{z^n}{(n+a)^s},\quad(|z|\leq1;\;a\neq0,-1,-2,\cdots),
\end{equation}
so that
\begin{equation}\label{yd2'}
\zeta(s,1)=\zeta(s),\; \zeta(s,2)=\zeta(s)-1 \quad\mbox{and}\quad
\zeta(s,a)+\Phi(-1,s,a)=\frac{1}{2^{s-1}}\zeta(s,\frac{a}{2}),
\end{equation}
are known to be meromorphic (that is, analytic everywhere in the
complex $s$-plane except for a simple pole at $s=1$ with residue 1).

There is also the multiplication theorem (see also reference
\cite{hurw})
\begin{equation}\label{yd1'}
m^s\zeta(s)=\sum^m_{k=1}\zeta(s,\frac{k}{m}),\quad(m\in\mathbb{N}),
\end{equation}
of which a useful generalization is
\begin{equation}\label{yd1''}
\sum^{m-1}_{k=0}\zeta(s,a+\frac{k}{m})=m^s\zeta(s,ma),\quad(m\in\mathbb{N}),
\end{equation}
so that
\begin{equation}\label{yd1t}
\zeta(s,1/3)+\zeta(s,2/3)=(3^s-1)\zeta(s)\quad\mbox{and}\quad\zeta(s,2/3)+\zeta(s,4/3)=3^s\zeta(s,2)-\zeta(s),
\end{equation}
are known to be meromorphic.

In Theorem 16, taking $x=c/3$, we can apply the identities
(\ref{yd2'}) and (\ref{yd1t}) in order to prove the following series
representations for $\zeta(2n+1)$:

\textbf{Theorem 17.}  Let $r\in\mathbb{N}$ be an arbitrary natural
number, and $B_k^*$ be the Bernoulli numbers, then the Riemann Zeta
function $\zeta(2n+1)$ can be recursively determined by the
following recurrence formula, namely

\parbox{11cm}{\begin{eqnarray*}\label{41}
 \zeta(2r+1)&=& \frac{2^{2r+1}3^{2r}}{3^{2r}(2^{2r}+1)+2^{2r}-1}\sum^{r-1}_{k=1}(-1)^{k-1}\frac{1}{(2k)!}\left(\frac{\pi}{3}\right)^{2k}\zeta(2r+1-2k)\\
 & &
 +\,(-1)^{r-1}\frac{2^{2r+1}\pi^{2r}}{(3^{2r}(2^{2r}+1)+2^{2r}-1)(2r)!}\left(H_{2r}-\ln\frac{\pi}{3}\right)\\
 & &
 +\,(-1)^{r-1}\frac{4(2\pi)^{2r}}{3^{2r}(2^{2r}+1)+2^{2r}-1}\sum^\infty_{k=1}\frac{(2k-1)!}{(2r+2k)!}\frac{\zeta(2k)}{6^{2k}},
\end{eqnarray*}}\hfill\parbox{1cm}{\begin{eqnarray}\end{eqnarray}}
or, equivalently,

\parbox{11cm}{\begin{eqnarray*}\label{42}
 \zeta(2r+1)&=& \frac{2^{2r+1}3^{2r}}{3^{2r}(2^{2r}+1)+2^{2r}-1}\sum^{r-1}_{k=1}(-1)^{k-1}\frac{1}{(2k)!}\left(\frac{\pi}{3}\right)^{2k}\zeta(2r+1-2k)\\
 & &
 +\,(-1)^{r-1}\frac{2^{2r+1}\pi^{2r}}{(3^{2r}(2^{2r}+1)+2^{2r}-1)(2r)!}\left(H_{2r}-\ln\frac{\pi}{3}\right)\\
 & &
 +\,(-1)^{r-1}\frac{(2\pi)^{2r}}{3^{2r}(2^{2r}+1)+2^{2r}-1}\sum^\infty_{k=1}\frac{B_k^*}{(2r+2k)!k}\left(\frac{\pi}{3}\right)^{2k}.
\end{eqnarray*}}\hfill\parbox{1cm}{\begin{eqnarray}\end{eqnarray}}
Obviously (\ref{41}) and (\ref{40''}) are identical.

\textbf{Proof.} Taking $x=c/3$ in Theorem 16, we obtain

\parbox{11cm}{\begin{eqnarray*}\label{41'}
 \sum^\infty_{n=1}\frac{1}{n^{2r+1}}\cos\frac{n\pi}{3}
 &=& \sum^{r-1}_{k=0}(-1)^k\frac{1}{(2k)!}\left(\frac{\pi}{3}\right)^{2k}\zeta(2r+1-2k)\\
 & &
 +\,(-1)^r\frac{1}{(2r)!}\left(\frac{\pi}{3}\right)^{2r}\left(H_{2r}-\ln\frac{\pi}{3}\right)\\
 & & +\,(-1)^r\left(\frac{\pi}{3}\right)^{2r}\sum^\infty_{k=1}\frac{(2k-1)!\,\zeta(2k)}{2^{2k-1}(2r+2k)!}\left(\frac{1}{3}\right)^{2k}.
\end{eqnarray*}}\hfill\parbox{1cm}{\begin{eqnarray}\end{eqnarray}}
Where
\begin{eqnarray*}
   \sum^\infty_{n=1}\frac{1}{n^{2r+1}}\cos\frac{n\pi}{3} &=& \cos\frac{\pi}{3}+ \sum^\infty_{n=1}\frac{1}{(3n-1)^{2r+1}}\cos\left(\frac{3n\pi}{3}-\frac{\pi}{3}\right)\\
   & & + \,\sum^\infty_{n=1}\frac{1}{(3n)^{2r+1}}\cos\frac{3n\pi}{3}+\sum^\infty_{n=1}\frac{1}{(3n+1)^{2r+1}}\cos\left(\frac{3n\pi}{3}+\frac{\pi}{3}\right) \\
   &=&
   \frac{1}{2}-\frac{\Phi(-1,2r+1,2/3)}{2\times3^{2r+1}}-\frac{\eta(2r+1)}{3^{2r+1}}-\frac{\Phi(-1,2r+1,4/3)}{2\times3^{2r+1}}\\
   &=&
   \frac{1}{2}-\frac{2^{-2r}\zeta(2r+1,1/3)-\zeta(2r+1,2/3)}{2\times3^{2r+1}}\\
   & &-\,\frac{(1-2^{-2r})\zeta(2r+1)}{3^{2r+1}}-\frac{2^{-2r}\zeta(2r+1,2/3)-\zeta(2r+1,4/3)}{2\times3^{2r+1}}\\
   &=& \frac{1}{2}-\frac{\zeta(2r+1,1/3)+\zeta(2r+1,2/3)}{2^{2r+1}3^{2r+1}}-\frac{(1-2^{-2r})\zeta(2r+1)}{3^{2r+1}}\\
   & &+\,\frac{\zeta(2r+1,2/3)+\zeta(2r+1,4/3)}{2\times3^{2r+1}}\\
   &=& \frac{1}{2}-\frac{(3^{2r+1}-1)\zeta(2r+1)}{2^{2r+1}3^{2r+1}}-\frac{(1-2^{-2r})\zeta(2r+1)}{3^{2r+1}}\\
   & &+\,\frac{3^{2r+1}[\zeta(2r+1)-1]-\zeta(2r+1)}{2\times3^{2r+1}}\\
   &=&
   \left(\frac{1}{2}+\frac{1}{2^{2r+1}3^{2r}}-\frac{1}{2^{2r+1}}-\frac{1}{2\times3^{2r}}\right)\zeta(2r+1).
\end{eqnarray*}
So the (\ref{41}) is proved.

Connon \cite{conn}, Srivastava and Tsumura \cite{sri} reported for
$\Re(s)>1$
\begin{equation}\label{yc1}
\sum^\infty_{n=1}\frac{1}{n^s}\cos\frac{n\pi}{3}=\frac{1}{2}(6^{1-s}-3^{1-s}-2^{1-s}+1)\zeta(s),
\end{equation}
\begin{equation}\label{yc2}
\sum^\infty_{n=1}\frac{1}{n^s}\cos\frac{2n\pi}{3}=\frac{1}{2}(3^{1-s}-1)\zeta(s),
\end{equation}
\begin{equation}\label{yc3}
\sum^\infty_{n=1}\frac{1}{n^s}\sin\frac{2n\pi}{3}=\sqrt{3}\left\{\frac{3^{-s}-1}{2}\zeta(s)+3^{-s}\zeta(s,\frac{1}{3})\right\},
\end{equation}
\begin{equation}\label{yc4}
\sum^\infty_{n=1}\frac{1}{n^s}\cos\frac{n\pi}{2}=2^{-s}(2^{1-s}-1)\zeta(s),
\end{equation}

Taking $x=2c/3$ in (\ref{33}), we can apply the (\ref{yc3}) in order
to prove the following Corollary:

\textbf{Corollary 3.} For Hurwitz Zeta function and Riemann Zeta
functions, we have following identities ($r\in\mathbb{N}$):

\parbox{11cm}{\begin{eqnarray*}\label{yc3'}
 & &
 \sqrt{3}\left[\zeta(2r+1,\frac{1}{3})+\frac{1-3^{2r+1}}{2}\,\zeta(2r+1)\right]\\
 &=&
 \sum^{r-1}_{k=0}(-1)^{k}\frac{(2\pi)^{2k+1}}{(2k+1)!}3^{2r-2k}\zeta(2r-2k)
 +\frac{(-1)^{r}2^{2r-1}\pi^{2r+1}}{(2r+1)!}(6r+1).
\end{eqnarray*}}\hfill\parbox{1cm}{\begin{eqnarray}\end{eqnarray}}

\textbf{Theorem 18.}  For $r\in\mathbb{N}$, in the interval $[-c,c]$
we have the following Fourier series expressions:

\parbox{11cm}{\begin{eqnarray*}\label{44}
 & &
 \sum^{\infty}_{n=1}(-1)^{n-1}\frac{1}{n^{2r}}\cos\frac{n\pi{x}}{c}\\
 &=&
 \sum^{r-1}_{k=0}(-1)^{k}\frac{1}{(2k)!}\left(\frac{\pi{x}}{c}\right)^{2k}\eta(2r-2k)
 +(-1)^{r}\frac{1}{2}\frac{1}{(2r)!}\left(\frac{\pi{x}}{c}\right)^{2r}.
\end{eqnarray*}}\hfill\parbox{1cm}{\begin{eqnarray}\end{eqnarray}}
\parbox{11cm}{\begin{eqnarray*}\label{45}
 & &
 \sum^{\infty}_{n=1}(-1)^{n-1}\frac{1}{n^{2r+1}}\sin\frac{n\pi{x}}{c}\\
 &=&
 \sum^{r}_{k=1}(-1)^{k-1}\frac{1}{(2k-1)!}\left(\frac{\pi{x}}{c}\right)^{2k-1}\eta(2r+2-2k)\\
 & & +\,(-1)^r\frac{1}{2}\frac{1}{(2r+1)!}\left(\frac{\pi{x}}{c}\right)^{2r+1}.
\end{eqnarray*}}\hfill\parbox{1cm}{\begin{eqnarray}\end{eqnarray}}

\textbf{Proof.}  According to Definition 1, if $S_0(t)=\ln(1+t)$,
then
\[S_{m-1}(t)=\sum^\infty_{n=1}(-1)^{n-1}\frac{t^n}{n^m},\quad{S_{m-1}(1)}=\eta(m).\]
In Theorem 7, let $m=2r$, then $S_{2r-1-2k}(1)=\eta(2r-2k)$. By
using Theorem 2 and Lemma 4, thus the (\ref{44}) is proved.

In Theorem 9, let $m=2r+1$, then $S_{2r+1-2k}(1)=\eta(2r+2-2k)$. By
using Theorem 2 and Lemma 4, thus the (\ref{45}) is proved.

\textbf{Theorem 19.}  The Dirichlet Eta function $\eta(2n)$
satisfies the following recurrence formula:
\begin{equation}\label{46}
    \sum^{r-1}_{k=0}(-1)^{k}\frac{\pi^{2k}}{(2k+1)!}\eta(2r-2k)=(-1)^{r-1}\frac{1}{2}\frac{\pi^{2r}}{(2r+1)!}.
\end{equation}

\textbf{Proof.}  Let $x=c$ in the second expression of Theorem 18,
then it can be proved.

\section{Sum-$a_K\pi^K$ of the Dirichlet series}
\noindent

We have already known many Dirichlet series of integer variables
have a sum similar to the one given by
\begin{equation}\label{euler}
\zeta(2n)=\frac{2^{2n-1}B_n^*}{(2n)!}\pi^{2n},
\end{equation}
therefore, we introduce the definition of the sum-$a_K\pi^K$:

\textbf{Definition 2.}  The sum-$a_K\pi^K$ of the Dirichlet series
of integer variables are numbers like $a_K\pi^K$, where $K$ is a
natural number, and $a_K$ is a rational number.

Of course we expect the Riemann Zeta function $\zeta(2r+1)$ has the
sum-$a_K\pi^K$, but to our disappointment, the calculation formula
of $\zeta(2r+1)$ contains the infinite power series of $\pi$, though
it converges very fast. Then we would like to ask: Whether or not
$\zeta(2r+1)$ has a sum-$a_K\pi^K$?

\textbf{Theorem 20.} If the Dirichlet series of integer variables
denoted by $f(m)$ has a sum-$a_K\pi^K$, then the corresponding
cosine series definitely converges to the polynomial function in a
certain interval.

\textbf{Proof.}  Denoting the sum of the Dirichlet series of integer
variables as $f(m)$ , namely
\begin{equation}\label{04}
f(m)=\sum^\infty_{n=1}\frac{a_n}{n^m},\quad m\in\mathbb{N}.
\end{equation}
Especially, we have
\[f(m)=\left\{\begin{array}{lr}\zeta(m),& a_n=1,\\\eta(m),
& a_n=(-1)^{n-1}.\end{array}\right.\] If the corresponding cosine
series of $f(m)$ has an analytic sum function in the interval
$a<0\leq{x}<{b}$, we can expand this analytic function into a power
series in the neighborhood of $x=0$, denoted by $P(\pi{x}/c)$,
namely
\[\sum^\infty_{n=1}\frac{a_n}{n^m}\cos\frac{n\pi x}{c}=P\left(\frac{\pi x}{c}\right),\quad a<0\leq x<b.\]
Differentiating both sides of the equation to obtain derivatives of
all orders less than $m$, and let $x=0$, we derive that the first
several terms of coefficients of $P(\pi{x}/c)$ are related to
$f(m-2k),\;k=0,1,\ldots,[m/2]-1$, namely
\[P\left(\frac{\pi{x}}{c}\right)=\sum^{[m/2]-1}_{k=0}(-1)^kf(m-2k)\frac{1}{(2k)!}\left(\frac{\pi{x}}{c}\right)^{2k}+R\left(\frac{\pi x}{c}\right).\]
If $x=c$ is in the interval $(a,b)$, then differentiating at $x=c$
we have
\begin{equation}\label{42}
\sum^{[m/2]-1}_{k=1}(-1)^kf(m-2k)\frac{\pi^{2k-1}}{(2k-1)!}+R'(\pi)=0.
\end{equation}

We would like to ask: Is $R'(\pi)$ in (\ref{42}) unique? Is there
another power series $R'(\pi)$ of $\pi$ making (\ref{42}) tenable?
As we know, if the analytic function $g(x)$ is analytic in a
neighborhood of $x=0$, then in a certain interval $a\leq x\leq b$ of
$x=0$, we have the unique power series expansion:
\[g(x)=a_0+a_1x+a_2x^2+\cdots,\quad a\leq x\leq b.\]
In other words, at any point $x=x_0$ in the interval
$a\leq{x}\leq{b}$, the power series of the function $g(x_0)$ is
unique. According to such a uniqueness of power series expansions,
the power series $R'(\pi)$ of $\pi$ in (\ref{42}) is definitely
unique. Therefore, if $R'(\pi)$ is an infinite power series of
$\pi$, then it definitely cannot be a polynomial of $\pi$, vice
versa.

If the Dirichlet series $f(m)$ has a sum-$a_K\pi^K$, according to
(\ref{42}), $R'(\pi)$ definitely has the form of the sum-$a_K\pi^K$.
But unless $R(\pi{x}/c)$ is a polynomial, $R'(\pi)$ cannot be the
sum-$a_K\pi^K$, therefore, $R(\pi{x}/c)$ is definitely a polynomial
and cannot be an infinite power series at the same time, thus
$P(\pi{x}/c)$ is definitely a polynomial. Then Theorem 20 is proved.

\textbf{Theorem 21.}  Let $r\geq1$, and $B_n^*$ be the Bernoulli
numbers, then in the interval $[-c,c]$, we have the following
Fourier series expansion related to $\eta(2n+1)$:

\parbox{11cm}{\begin{eqnarray*}\label{43}
 & &
 \sum^\infty_{n=1}(-1)^{n-1}\frac{1}{n^{2r+1}}\cos\frac{n\pi x}{c}\\
 &=&
 \sum^{r}_{k=0}(-1)^{k}\frac{1}{(2k)!}\left(\frac{\pi x}{c}\right)^{2k}\eta(2r+1-2k)\\
 & & +\,(-1)^{r+1}\sum^\infty_{n=1}\frac{(2^{2n}-1)B_n^*}{2n(2r+2n)!}\left(\frac{\pi x}{c}\right)^{2r+2n}.
\end{eqnarray*}}\hfill\parbox{1cm}{\begin{eqnarray}\end{eqnarray}}

\textbf{Proof.}  According to Definition 1, if $S_0(t)=\ln(1+t)$,
then
\[S_{m-1}(t)=\sum^\infty_{n=1}(-1)^{n-1}\frac{t^{n}}{n^m},\quad S_{m-1}(1)=\eta(m).\]
In Theorem 10, let $m=2r+1$, using Theorem 2, and considering
$S_{2r-2k}(1)=\eta(2r+1-2k)$ we have
\begin{eqnarray*}
   & & \sum^\infty_{n=1}(-1)^{n-1}\frac{1}{n^{2r+1}}\cos\frac{n\pi x}{c} \\
   &=& \sum^{r-1}_{k=0}(-1)^{k}\frac{1}{(2k)!}\left(\frac{\pi x}{c}\right)^{2k}\eta(2r+1-2k) \\
   & &
   +\,(-1)^{r}\left(\frac{\pi}{c}\right)^{2r}\underbrace{\int^x_0dx\cdots}_{2r}\int^x_0
   \left.\cos\left(\frac{\pi x}{c}\frac{\partial}{\partial z}\right)\ln(1+e^z)\right|_{z=0}dx.
\end{eqnarray*}
In the above expression, by using the algorithms and basic formulas
in this paper we have
\begin{eqnarray*}
   & & \left.\cos\left(\frac{\pi x}{c}\frac{\partial}{\partial z}\right)\ln(1+e^z)\right|_{z=0}=\left.\frac{1}{2}\ln(X^2+Y^2)\right|_{z=0} \\
   &=&
   \frac{1}{2}\ln\left(\left(1+\cos\frac{\pi x}{c}\right)^2+\sin^2\frac{\pi x}{c}\right)=\ln\left(2\cos\frac{\pi x}{2c}\right)\\
   &=& \ln2-\sum^\infty_{n=1}\frac{2^{2n}(2^{2n}-1)B_n^*}{2n(2n)!}\left(\frac{\pi x}{2c}\right)^{2n},\quad|x|<c.
\end{eqnarray*}
Substituting this result into the above expression, and considering
$\eta(1)=\ln2$, then when $r\geq1$, the formula is tenable at the
endpoint $|x|=c$, thus we have (\ref{43}) immediately.

\textbf{Theorem 22.}  The Riemann Zeta function $\zeta(2n+1)$ of odd
variables does not have a sum-$a_K\pi^K$.

\textbf{Proof.} According to Theorem 20, the necessary condition for
the Dirichlet Eta function $\eta(2n+1)$ having the sum-$a_K\pi^K$ is
that the corresponding cosine series takes a polynomial as its sum
function in a certain interval. According to Theorem 21 and the
uniqueness of power series expansions, the corresponding cosine
series of $\eta(2n+1)$ has an analytic sum function only in the
interval $[-c,c]$, but such a sum function cannot be a polynomial,
thus the Dirichlet Eta function $\eta(2n+1)$ cannot have a
sum-$a_K\pi^K$. In addition,
\begin{equation}\label{05}
\zeta(2n+1)=\frac{2^{2n}}{2^{2n}-1}\eta(2n+1),\quad{n\in\mathbb{N}}.
\end{equation}
Therefore, the Riemann Zeta function $\zeta(2n+1)$ of odd variables
does not have a sum-$a_K\pi^K$ as well.

In precisely the same manner, we can prove the following results for
the Dirichlet Beta function $\beta(2n+1)$:

\textbf{Corollary 5.} For $r\in\mathbb{N}$, in the interval
$[-c/2,c/2]$ we have the following Fourier series expressions:

\parbox{11cm}{\begin{eqnarray*}\label{beta}
 & &
 \sum^{\infty}_{n=0}(-1)^n\frac{1}{(2n+1)^{2r+1}}\cos\frac{(2n+1)\pi{x}}{c}\\
 &=& \sum^{r-1}_{k=0}(-1)^{k}\frac{1}{(2k)!}\left(\frac{\pi{x}}{c}\right)^{2k}\beta(2r+1-2k)
+(-1)^r\frac{\pi}{4}\frac{1}{(2r)!}\left(\frac{\pi{x}}{c}\right)^{2r}.
\end{eqnarray*}}\hfill\parbox{1cm}{\begin{eqnarray}\end{eqnarray}}
Where the Dirichlet Beta function is defined as
\begin{equation}\label{yd3}
\beta(s)=\sum^\infty_{n=0}\frac{(-1)^n}{(2n+1)^s},\quad(\Re(s)>0).
\end{equation}

By using the (\ref{beta}) we can easily obtain
\begin{equation}\label{beta'}
\sum^{r-1}_{k=0}(-1)^{k}\frac{1}{(2k)!}\left(\frac{\pi}{2}\right)^{2k}\beta(2r+1-2k)
=(-1)^{r-1}\frac{\pi}{4}\frac{1}{(2r)!}\left(\frac{\pi}{2}\right)^{2r},\quad(r\in\mathbb{N}).
\end{equation}

For any positive integer k:
\begin{equation}\label{ybeta}
\beta(2k+1)=\frac{(-1)^kE_{2k}\pi^{2k+1}}{4^{k+1}(2k)!},
\end{equation}
where $E_n$ represent the Euler numbers. By using the (\ref{beta'})
and (\ref{ybeta}) we can easily obtain
\begin{equation}\label{ybeta'}
\sum^{r-1}_{k=0}{2r\choose{2k}}E_{2r-2k}=-1,\quad{r\in\mathbb{N}}.
\end{equation}

According to Definition 1, if $S_0(t)=\ln\frac{1+t}{1-t}$, then
\[S_{m-1}(t)=\sum^\infty_{n=1}\frac{t^{2n-1}}{(2n-1)^m},\quad{S_{m-1}(1)}=\lambda(m).\]
In Theorem 7, let $m=2r$, and by using Theorem 2, considering
$S_{2r-1-2k}(1)=\lambda(2r-2k)$ and
\[\left.\sin\left(\frac{\pi{x}}{c}\frac{\partial}{\partial{z}}\right)\left(\frac{1}{2}\ln\frac{1+e^z}{1-e^z}\right)\right|_{z=0}=
 \frac{\pi}{4},\quad0<x<c,\]
thus we obtain

\textbf{Corollary 6.} For $r\in\mathbb{N}$, in the interval $[0,c]$
we have

\parbox{11cm}{\begin{eqnarray*}\label{044}
 & &
 \sum^{\infty}_{n=1}\frac{1}{(2n-1)^{2r}}\cos\frac{(2n-1)\pi{x}}{c}\\
 &=& \sum^{r-1}_{k=0}(-1)^{k}\frac{1}{(2k)!}\left(\frac{\pi{x}}{c}\right)^{2k}\lambda(2r-2k)
+(-1)^r\frac{\pi}{4}\frac{1}{(2r-1)!}\left(\frac{\pi{x}}{c}\right)^{2r-1}.
\end{eqnarray*}}\hfill\parbox{1cm}{\begin{eqnarray}\end{eqnarray}}
Where $\lambda(s)$ are the Dirichlet Lambda function defined by
\begin{equation}\label{yd5}
\lambda(s)=\sum^\infty_{n=0}\frac{1}{(2n+1)^s},\quad(\Re(s)>1).
\end{equation}
Letting $x=c/2$, we obtain
\begin{equation}\label{044'}
\sum^{r-1}_{k=0}(-1)^k\frac{1}{(2k)!}\left(\frac{\pi}{2}\right)^{2k}\lambda(2r-2k)=\frac{(-1)^{r-1}}{2(2r-1)!}\left(\frac{\pi}{2}\right)^{2r},\quad{r\in\mathbb{N}}.
\end{equation}

Letting $x=c/4$ in (\ref{044}), we obtain

\parbox{11cm}{\begin{eqnarray*}\label{56}
 & &
 \sum^\infty_{n=1}(-1)^{[\frac{n}{2}]}\frac{1}{(2n-1)^{2r}}\\
 &=&
 \sum^{r-1}_{k=0}(-1)^k\frac{\sqrt{2}}{(2k)!}\left(\frac{\pi}{4}\right)^{2k}\lambda(2r-2k)
 +(-1)^r\frac{\sqrt{2}}{(2r-1)!}\left(\frac{\pi}{4}\right)^{2r}.
\end{eqnarray*}}\hfill\parbox{1cm}{\begin{eqnarray}\end{eqnarray}}
For example, letting $r=1,2$ in (\ref{56}), and using (\ref{044'})
we obtain

\parbox{11cm}{\begin{eqnarray*}\label{044''}
\sum^\infty_{n=1}(-1)^{[\frac{n}{2}]}\frac{1}{(2n-1)^2}&=&\frac{\pi^2}{16}\sqrt{2},\\
\sum^\infty_{n=1}(-1)^{[\frac{n}{2}]}\frac{1}{(2n-1)^4}&=&\frac{11\pi^4}{1536}\sqrt{2}.
\end{eqnarray*}}\hfill\parbox{1cm}{\begin{eqnarray}\end{eqnarray}}

\textbf{Theorem 23} Let
$S(t)=\sum^\infty_{n=0}a_nt^n,\;t\in\mathbb{R}^1,\;0\leq{t}\leq{r},\;0<r<+\infty$,
if
\[\sum^\infty_{n=0}a_n\cos\frac{n\pi{x}}{c}=\left.\cos\left(\frac{\pi{x}}{c}\frac{\partial}{\partial{z}}\right)S(e^z)\right|_{z=0},\]
\[\sum^\infty_{n=0}a_n\sin\frac{n\pi{x}}{c}=\left.\sin\left(\frac{\pi{x}}{c}\frac{\partial}{\partial{z}}\right)S(e^z)\right|_{z=0}\]
are tenable in $a<x<b,\;x\in\mathbb{R}^1$, then
$\forall{x_0}\in[0,\frac{b-a}{2})$,

\parbox{11cm}{\begin{eqnarray*}\label{58}
 \sum^\infty_{n=0}a_n\cos\frac{n\pi{x_0}}{c}\cos\frac{n\pi{x}}{c}&=&
 \cosh\left(x_0\frac{\partial}{\partial{x}}\right)\left.\cos\left(\frac{\pi{x}}{c}\frac{\partial}{\partial{z}}\right)S(e^z)\right|_{z=0},\\
 \sum^\infty_{n=0}a_n\cos\frac{n\pi{x_0}}{c}\sin\frac{n\pi{x}}{c}&=&
 \cosh\left(x_0\frac{\partial}{\partial{x}}\right)\left.\sin\left(\frac{\pi{x}}{c}\frac{\partial}{\partial{z}}\right)S(e^z)\right|_{z=0}
\end{eqnarray*}}\hfill\parbox{1cm}{\begin{eqnarray}\end{eqnarray}}
are tenable in $a+x_0<x<b-x_0$. Accordingly, for any definite value
$x\in(a,b)$, $x_0$ in (\ref{58}) takes values in the following
interval:
\[0<x_0<\left\{\begin{array}{lr}x-a,&a<x\leq(a+b)/2,\\b-x,&b>x\geq(a+b)/2.\end{array}\right.\]

\textbf{Proof}  Clearly
\begin{eqnarray*}
   & &
\sum^\infty_{n=0}a_n\cos\frac{n\pi{x_0}}{c}\cos\frac{n\pi{x}}{c}=
\left.\cos\left(\frac{\pi{x_0}}{c}\frac{\partial}{\partial{z}}\right)\cos\left(\frac{\pi{x}}{c}\frac{\partial}{\partial{z}}\right)S(e^z)\right|_{z=0} \\
   &=& \frac{1}{2}\left.\cos\left(\frac{\pi(x-x_0)}{c}\frac{\partial}{\partial{z}}\right)S(e^z)\right|_{z=0}
   +\frac{1}{2}\left.\cos\left(\frac{\pi(x+x_0)}{c}\frac{\partial}{\partial{z}}\right)S(e^z)\right|_{z=0} \\
   &=& \cosh\left(x_0\frac{\partial}{\partial{x}}\right)\left.\cos\left(\frac{\pi{x}}{c}\frac{\partial}{\partial{z}}\right)S(e^z)\right|_{z=0}.
\end{eqnarray*}
\begin{eqnarray*}
   & &
\sum^\infty_{n=0}a_n\cos\frac{n\pi{x_0}}{c}\sin\frac{n\pi{x}}{c}=
\left.\cos\left(\frac{\pi{x_0}}{c}\frac{\partial}{\partial{z}}\right)\sin\left(\frac{\pi{x}}{c}\frac{\partial}{\partial{z}}\right)S(e^z)\right|_{z=0} \\
   &=& \frac{1}{2}\left.\sin\left(\frac{\pi(x-x_0)}{c}\frac{\partial}{\partial{z}}\right)S(e^z)\right|_{z=0}
   +\frac{1}{2}\left.\sin\left(\frac{\pi(x+x_0)}{c}\frac{\partial}{\partial{z}}\right)S(e^z)\right|_{z=0} \\
   &=& \cosh\left(x_0\frac{\partial}{\partial{x}}\right)\left.\sin\left(\frac{\pi{x}}{c}\frac{\partial}{\partial{z}}\right)S(e^z)\right|_{z=0}.
\end{eqnarray*}
These two equations are both tenable in the common interval
$a+x_0<x<b-x_0$ of $a<x-x_0<b$ and $a<x+x_0<b$.

For (\ref{044}), using Theorem 23 we can easily obtain:

For any definite value $x_0\in[0,c/2)$, we have

\parbox{12cm}{\begin{eqnarray*}\label{59}
   & & \sum^\infty_{n=1}\frac{1}{(2n-1)^{2r}}\cos\frac{(2n-1)\pi{x_0}}{c}\cos\frac{(2n-1)\pi{x}}{c} \\
   &=& \sum^{r-1}_{k=0}(-1)^k\frac{1}{(2k)!}\,\lambda(2r-2k)\,\cosh\left(x_0\frac{\partial}{\partial{x}}\right)\left(\frac{\pi{x}}{c}\right)^{2k} \\
   & &
   +(-1)^r\frac{1}{(2r-1)!}\frac{\pi}{4}\,\cosh\left(x_0\frac{\partial}{\partial{x}}\right)\left(\frac{\pi{x}}{c}\right)^{2r-1}
\end{eqnarray*}}\hfill\parbox{1cm}{\begin{eqnarray}\end{eqnarray}}
are tenable in $x_0\leq{x}\leq{c-x_0}$. On the contrary, for any
definite value $x\in(0,c)$, $x_0$ in (\ref{59}) takes value in the
following interval:
\[0\leq{x_0}\leq\left\{\begin{array}{lr}x,&0<x\leq{c/2},\\c-x,&c>x\geq{c/2}.\end{array}\right.\]

When $r=1$, the equation above can be expressed as:
\[\sum^\infty_{n=1}\frac{1}{(2n-1)^{2}}\cos\frac{(2n-1)\pi{x_0}}{c}\cos\frac{(2n-1)\pi{x}}{c}=\lambda(2)-\frac{\pi^2x}{4c}.\]

As the equation is tenable in $x_0\leq{x}\leq{c-x_0}$, $x$ cannot be
$0$ and $c$ unless $x_0=0$. When $x_0=c/4\in[0,c/2)$, we get
\[\sum^\infty_{n=1}(-1)^{[n/2]}\frac{1}{(2n-1)^{2}}\cos\frac{(2n-1)\pi{x}}{c}=\lambda(2)\sqrt{2}-\frac{\pi^2x}{4c}\sqrt{2}\]
is tenable in $c/4\leq{x}\leq{3c/4}$, where $\lambda(2)=\pi^2/8$.

If we take not $x_0=c/4$ but $x=c/4\in(0,c/2]$, then in
$0\leq{x_0}\leq{c/4}$ we have
\[\sum^\infty_{n=1}(-1)^{[n/2]}\frac{1}{(2n-1)^{2}}\cos\frac{(2n-1)\pi{x_0}}{c}=\lambda(2)\sqrt{2}-\frac{\pi^2}{16}\sqrt{2}=\frac{\pi^2}{16}\sqrt{2}.\]
Combining these two equations into one, we have
\[\sum^\infty_{n=1}(-1)^{[\frac{n}{2}]}\frac{1}{(2n-1)^{2}}\cos\frac{(2n-1)\pi{x}}{c}
=\left\{\begin{array}{lr}\pi^2\sqrt{2}/16,&0\leq{x}\leq{c}/4,\\\lambda(2)\sqrt{2}-\pi^2x\sqrt{2}/(4c),
&c/4\leq{x}\leq3c/4.\end{array}\right.\]

In (\ref{59}), if we take not $r=1$ but $r=2$, then similarly we
have

\parbox{11cm}{\begin{eqnarray*}\label{69}
 & &
 \frac{1}{\sqrt{2}}\sum^\infty_{n=1}(-1)^{[n/2]}\frac{1}{(2n-1)^{4}}\cos\frac{(2n-1)\pi x}{c}\\
 &=& \frac{5\pi^4}{768}+\frac{\pi^4x}{128c}-\frac{\pi^4x^2}{16c^2}+\frac{\pi^4x^3}{24c^3},\quad\frac{c}{4}\leq x\leq\frac{3c}{4}.
\end{eqnarray*}}\hfill\parbox{1cm}{\begin{eqnarray}\end{eqnarray}}

\parbox{11cm}{\begin{eqnarray*}\label{70}
 & &
 \frac{1}{\sqrt{2}}\sum^\infty_{n=1}(-1)^{[n/2]}\frac{1}{(2n-1)^{4}}\cos\frac{(2n-1)\pi x_0}{c}\\
 &=& \frac{11\pi^4}{1536}-\frac{\pi^4}{32c^2}x_0^2,\quad0\leq x_0\leq{c}/4.
\end{eqnarray*}}\hfill\parbox{1cm}{\begin{eqnarray}\end{eqnarray}}

If the following Dirichlet series is denoted by $\mathfrak{D}(2r)$,
\begin{equation}\label{yd6}
\mathfrak{D}(2r)=\frac{1}{\sqrt{2}}\sum^\infty_{n=1}(-1)^{[n/2]}\frac{1}{(2n-1)^{2r}},\quad{r\in\mathbb{N}},
\end{equation}
then by using (\ref{044}) and Theorem 20, we know that
$\mathfrak{D}(2r)$ has a sum-$a_K\pi^K$. The (\ref{044''}) gives
$\mathfrak{D}(4)=11\pi^4/1536,\;\mathfrak{D}(2)=\pi^2/16$. The
(\ref{69}) indicates that the corresponding cosine series of
$\mathfrak{D}(4)$ converge to polynomials. By observing (\ref{70})
we can find that coefficients of the polynomial functions on the
right side of the equation are definitely related to
$\mathfrak{D}(4)$ and $\mathfrak{D}(2)$, and can be expressed as:
\[\frac{1}{\sqrt{2}}\sum^\infty_{n=1}(-1)^{[n/2]}\frac{1}{(2n-1)^{4}}\cos\frac{(2n-1)\pi x_0}{c}
=\mathfrak{D}(4)-\frac{1}{2!}\mathfrak{D}(2)\left(\frac{\pi
x_0}{c}\right)^2.\] The result is the same as that pointed out in
the proof of Theorem 20.

Generally, we can easily prove:

\textbf{Corollary 7.} For $r\in\mathbb{N}$, in the interval
$[0,c/4]$ we have

\parbox{11cm}{\begin{eqnarray*}\label{71}
 & &
 \frac{1}{\sqrt{2}}\sum^\infty_{n=1}(-1)^{[n/2]}\frac{1}{(2n-1)^{2r}}\cos\frac{(2n-1)\pi x_0}{c}\\
 &=& \sum^{r-1}_{k=0}(-1)^k\frac{1}{(2k)!}\left(\frac{\pi x_0}{c}\right)^{2k}\mathfrak{D}(2r-2k),\quad0\leq x_0\leq{c}/4.
\end{eqnarray*}}\hfill\parbox{1cm}{\begin{eqnarray}\end{eqnarray}}

Letting $x_0=c/4$ in (\ref{71}), we get the recurrence formula of
$\mathfrak{D}(2r)$ as

\begin{equation}\label{73}
\sum^{r-1}_{k=0}(-1)^k\frac{1}{(2k)!}\left(\frac{\pi}{4}\right)^{2k}\mathfrak{D}(2r-2k)=\frac{1}{2}\lambda(2r).
\end{equation}
Clearly $\lambda(2r)$ satisfies
$\zeta(2r)=\lambda(2r)+\zeta(2r)/2^{2r}$, then by solving this
equation for $\lambda(2r)$ and substituting it into (\ref{73}), we
have

\begin{equation}\label{72}
\sum^{r-1}_{k=0}(-1)^k\frac{1}{(2k)!}\left(\frac{\pi}{4}\right)^{2k}\mathfrak{D}(2r-2k)=\frac{2^{2r}-1}{2^{2r+1}}\zeta(2r).
\end{equation}

Using (\ref{72}) or (\ref{73}), we can recursively obtain the
sum-$a_K\pi^K$ of $\mathfrak{D}(2r)$. For instance, let $r=3$ in
(\ref{72}), using the sum-$a_K\pi^K$ of
$\mathfrak{D}(2),\;\mathfrak{D}(4)$ and
 $\zeta(6)$, we can easily obtain the
sum-$a_K\pi^K$ of $\mathfrak{D}(6)$:
\begin{equation}\label{74}
\mathfrak{D}(6)=\frac{1}{\sqrt{2}}\sum^\infty_{n=1}(-1)^{[n/2]}\frac{1}{(2n-1)^6}=\frac{361\pi^6}{491520}.
\end{equation}
The result is the same as that obtained by (\ref{56}).

In precisely the same manner, we can prove the following results:

\textbf{Corollary 8.} If the following Dirichlet series is denoted
by $\mathcal {D}(2r+1)$, namely
\begin{equation}\label{0yd6}
\mathcal{D}(2r+1)=\frac{1}{\sqrt{2}}\sum^\infty_{n=0}(-1)^{[n/2]}\frac{1}{(2n+1)^{2r+1}},\quad{r\in\mathbb{N}_0},
\end{equation}
then for $r\in\mathbb{N}$, in the interval $[0,c/4]$ we have

\parbox{11cm}{\begin{eqnarray*}\label{071}
 & &
 \frac{1}{\sqrt{2}}\sum^\infty_{n=0}(-1)^{[n/2]}\frac{1}{(2n+1)^{2r+1}}\cos\frac{(2n+1)\pi x_0}{c}\\
 &=& \sum^r_{k=0}(-1)^k\frac{1}{(2k)!}\left(\frac{\pi x_0}{c}\right)^{2k}\mathcal {D}(2r+1-2k),\quad0\leq x_0\leq{c}/4.
\end{eqnarray*}}\hfill\parbox{1cm}{\begin{eqnarray}\end{eqnarray}}
Where $\mathcal {D}(1)=\pi/4$, and $\forall{r}\in\mathbb{N}$
\begin{equation}\label{071'}
\mathcal{D}(2r+1)=\sum^{r-1}_{k=0}\frac{(-1)^k}{(2k+1)!}\left(\frac{\pi}{4}\right)^{2k+1}\lambda(2r-2k)+\frac{(-1)^r}{(2r)!}\left(\frac{\pi}{4}\right)^{2r+1}.
\end{equation}

Letting $x_0=c/4$ in (\ref{071}), we get the recurrence formula of
$\mathcal{D}(2r+1)$ as
\begin{equation}\label{071''}
\frac{1}{2}\beta(2r+1)=\sum^r_{k=0}\frac{(-1)^k}{(2k)!}\left(\frac{\pi}{4}\right)^{2k}\mathcal{D}(2r+1-2k),\quad{r}\in\mathbb{N}.
\end{equation}

For example, letting $r=1,2,3$ in (\ref{071''}), and using
(\ref{beta'}) we obtain
\begin{equation}\label{0071}
\mathcal{D}(3)=\frac{3\pi^3}{128},\quad\mathcal{D}(5)=\frac{57\pi^5}{24576},\quad\mathcal{D}(7)=\frac{307\pi^7}{1310720}.
\end{equation}
The results is the same as that obtained by (\ref{071'}).

Similarly, differentiating (\ref{071}) at $x_0=c/4$ we get
\begin{equation}\label{w71'}
\frac{1}{2}\lambda(2r)=\sum^{r-1}_{k=0}\frac{(-1)^k}{(2k+1)!}\left(\frac{\pi}{4}\right)^{2k+1}\mathcal{D}(2r-1-2k),\quad{r}\in\mathbb{N}.
\end{equation}

\textbf{Acknowledgments:} We are grateful to Wenpeng Zhang for some
helpful remarks. The authors also would like to take this
opportunity to thank Armen Bagdasaryan for his generous help.

\end{CJK*}
\end{document}